\documentclass[twocolumn,sts,preprint]{imsart}

\usepackage{amsthm,amsmath,amssymb}
\usepackage{amscd}
\usepackage{natbib}
\usepackage{relsize}
\usepackage{lineno}
\usepackage{etex}
\usepackage{pictex}

\startlocaldefs

\newcommand{\bi}[1]{\mathbf{#1}}
\newcommand{\bs}[1]{\boldsymbol{#1}}
\newtheorem{theo}{Proposition}
\newtheorem{rul}{Rule}
\newtheorem{lemm}{Lemma}
\newtheorem{define}{Definition}
\newtheorem{exam}{Example}
\newtheorem{remark}{Remark}
\newtheorem{proper}{Property}
\endlocaldefs

\begin{document}

\begin{frontmatter}

\title{On Probabilistic Parametric Inference}
\runtitle{Probabilistic Parametric Inference}

\author{\fnms{Toma\v z}
  \snm{Podobnik$^{1,2}$}\ead[label=e1]{Tomaz.Podobnik@ijs.si}}
\and
\author{\fnms{Tomi} \snm{\v Zivko$^2$}\ead[label=e2]{Tomi.Zivko@ijs.si}}

\runauthor{T. Podobnik and T. \v Zivko}

\affiliation{$^1$Faculty of Mathematics and Physics, University of Ljubljana 
      and $^2$Jo\v zef Stefan Institute, Ljubljana, Slovenia}

\address{\printead{e1,e2}}

\begin{abstract}
      An objective operational theory of probabilistic parametric inference
      is formulated without invoking the so-called non-informative
      prior probability distributions.
\end{abstract}



\end{frontmatter}

\section{Introduction}
\label{sec:introduction}
We make a probabilistic inference about a parameter of a
family of the so-called direct probability 
distributions by specifying a probability distribution that 
corresponds to the distribution of our belief in different values 
of the parameter (\cite{jef1}, \S\,2.0, p.\,22). 
The probabilistic parametric inference is
characteristic of {\em Bayesian schools} of statistical inference (as
opposed to {\em frequentist schools}), where the name Bayesian is
due to the central role of Bayes' Theorem in the process of inference. 
In the Bayesian paradigms, it is also possible to make 
statements concerning the values of the inferred parameters
in the absence of data, and these statements can be summarized 
in the so-called ({\em non-informative}) 
{\em prior probability distributions}, (\cite{vil}; see also, 
for example, 
\cite{jef}, \S\,1.4, p.\,33 and \S\,3.1,
pp.\,117-118; 
\cite{fer}, \S\,1.6, pp.\,30-31; 
\cite{berg}, \S\,1.2, pp.\,4-5;
\cite{rao}, \S\,3.5, p.\,86;
\cite{hag}, \S\,1.21, p,\,23;
\cite{kas};
\cite{lad}, \S\,3.4, p.\,150;
\cite{sha}, \S\,4.1.1, p.\,193;
\cite{rob}, \S\,3.5, pp.\,127-140;
\cite{case}, \S\,7.2.3, p.\,324;
\cite{jay}, \S\,4.1, pp.\,87-88;
\cite{har}, \S\,2.1, p.\,9;
\cite{hog}, \S\,11.2.1, pp.\,583-584).
The non-informative prior distributions provide a formal way of
expressing ignorance about the inferred parameter 
(\cite{jef}, \S\,3.1, pp.\,117-118; \cite{kas},
\S\,4.1, p.\,1355). It has been asserted  (\cite{jef1}, \S\,2.3,
p.\,31; \cite{jef}, \S\,1.5, pp.\,36-37; \cite{ber}, \S\,5.1, 
p.\,123; \cite{kas}, \S\,4.1, pp.\,1355-1356; 
\cite{rob}, \S\,3.5, p.\,127) that
there is no objective, unique non-informative prior distribution that 
represents ignorance. Instead, the priors should be chosen by public
agreement, much like units of length and weight, upon which everyone
could fall back when the prior information about the inferred
parameter is missing.

In the present article, a theory of probabilistic parametric 
inference is developed without invoking the non-informative prior probability 
distributions. Moreover, it is demonstrated that the non-informative 
prior probability distributions necessarily lead to inconsistencies.
Sections \ref{sec:direct}--\ref{sec:factors} are devoted to
formulation of a mathematical theory of probabilistic parametric
inference. In particular, in Section\;\ref{sec:direct}, the notions of 
probability, of (direct) probability distribution, of parametric
family and of invariant
family are introduced. In addition, some of the properties
of probability distributions are briefly reviewed. In
Section\;\ref{sec:inverse}, the so-called inverse probability
distributions are defined. It is demonstrated that the inverse
probability distributions must be directly proportional to the
appropriate direct probability distributions. The proportionality
factors, called consistency factors, are determined in 
Section\;\ref{sec:factors} on the grounds of
invariance of parametric families of direct probability distributions 
under the action of Lie groups. 
In Section\;\ref{sec:interpretations}, the
concept of relative frequency and the concept of degree of belief
are introduced that link the probability distributions to
an external world of measurable phenomena. In this way, the mathematical
theory becomes operational. Also in
Section\;\ref{sec:interpretations}, as well as in Conclusions, 
a reconciliation between
the Bayesian and the frequentist schools of parametric inference 
is advocated.
%
%

\section{Probabilities and probability distributions}
\label{sec:direct}

\subsection{Notation and general definitions}
\label{ss:general} 

In this section, the notions of probability and of probability
distribution are introduced, and some of the properties of
probability distributions are briefly reviewed, with special 
attention being paid to conditional probability density functions.
The purpose
of refreshing these well known concepts is to avoid misunderstandings
in subsequent sections where the properties of probability
distributions are extensively invoked and the definition of 
of the conditional probability distribution is extended.

Let $\Omega$ be a non-empty \textit{universal set}, also called
a  {\em sample space}, whose elements are
denoted by $\omega$. A set $\Sigma$ of subsets
$A,B,C\ldots$ of the sample space is called a
\textit{$\sigma$-algebra} (or \textit{$\sigma$-field}) 
on $\Omega$ if $\Sigma$ has $\Omega$ as a member, and is closed under
complementation, $\overline{A}\in \Sigma$; $\forall\,A\in\Sigma$, 
and under countable union,
$\sum_{i=1}^{\infty}A_i \in \Sigma$; $\forall\,A_1,A_2,\ldots\in\Sigma$ 
(throughout the present discussion,  $A+B$, $AB$ and $A-B$ denote a union,
an intersection and a relative complement of sets $A$ and $B$, 
 respectively, while $\overline{A}\equiv \Omega - A$). An ordered
pair $(\Omega,\Sigma)$ consisting of a state space $\Omega$ and 
a $\sigma$-algebra $\Sigma$ on $\Omega$ is called
a \textit{measurable space}. 
\begin{exam}[Borel algebra]
\label{ex:borel}
{\rm Let $\Omega$ be $\mathbb{R}^{n}$. The
  \textit{Borel $\sigma$-algebra} (or \textit{Borel algebra})
  ${\mathcal B}^n$ on
  $\mathbb{R}^{n}$ is the minimal $\sigma$-algebra containing a collection 
  of open rectangles in $\mathbb{R}^{n}$. It is also said 
  that the Borel algebra ${\mathcal B}^n$ on $\mathbb{R}^n$ is
 \textit{generated} by all open rectangles in $\mathbb{R}^{n}$.
 Every set from a Borel algebra is called a \textit{Borel set}.}
\end{exam}
\begin{define}[Probability]
\label{def:prob}
Let $P$ be a
real-valued function
on a $\sigma$-field $\Sigma$ on a sample space $\Omega$.
We call $P$ a  {\em probability measure} 
(or simply a {\em probability}) 
if it is congruent with the following three axioms due to
\cite{kol}: 
\begin{linenomath}
\begin{eqnarray}
 \label{eq:axiom1}
 P(A) & \ge & 0  \ ; \ \ \forall\ A\in\Sigma \ , \\
 \label{eq:axiom2} 
 P(\Omega) & = & 1 \ , \\
 \label{eq:axiom3}
 P\left(\sum_{i=1}^{\infty}A_i\right) & = & \sum_{i=1}^{\infty}P(A_i) 
\end{eqnarray}
\end{linenomath}
for all $A_i,A_j\in\Sigma$ that are \textit{mutually
  exclusive}, i.e., $A_iA_{j\neq i}=\emptyset$. Then, the triple 
$(\Omega,\Sigma, P)$ is termed the {\em probability space}.
\end{define}
\begin{define}[Random variable]
\label{def:random}
Given a probability space $(\Omega,\Sigma,P)$, let 
a function $X: \Omega \longrightarrow \mathbb{R}$ be
$\Sigma$-{\em measurable}:
$A_{X\le x}=\{\omega\in \Omega :X(\omega)\le x\} \in
\Sigma$, $\forall\,x\in\mathbb{R}$.
Then, $X$ is called a (real-valued) {\em scalar random variable} (or
{\em random variate}), while $x$ is called a 
{\em realization of} $X$.
\end{define}
\begin{define}[Distribution function]
\label{def:cdf}
Given a random variable $X$ on a probability space
$(\Omega,\Sigma,P)$, the ({\em cumulative}) {\em distribution 
function} (cdf) $F_X(x)$ is a real-valued function on 
the state space $\mathbb{R}$ to $[0,1]$ such that
$ F_X(x)=P(A_{X\le x})$.
\end{define}
Every cdf is a non-decreasing function with
$F_X(-\infty)\equiv\lim_{x\to -\infty}F_X(x)=0$ 
and
\begin{linenomath}
\begin{equation}
F_X(+\infty)\equiv\lim_{x\to +\infty}F_X(x)=1 \ . 
\label{eq:norcdf}
\end{equation}
\end{linenomath}
\begin{define}[Continuous random variable]
\label{def:pdf}
A random variable $X$ is called 
    {\em continuous} if its cdf $F_X(x)$ is {\em absolutely continuous},
    i.e., if the cdf is expressible as an integral of a non-negative
    (Lebesgue) integrable function $f_X(x)$, called 
    {\em probability density function} (pdf):
\begin{linenomath}
\begin{equation*}
 F_X(x)=\int_{-\infty}^{x}f_X(x')\,dx' \ .
\end{equation*}
\end{linenomath}
The {\em support} of a
continuous random variable $X$ is a set, say $V_X$, 
of all $x$ for which $f_X(x)>0$. 
\end{define}
Due to \eqref{eq:norcdf}, a pdf is always normalized to unit
area,
\begin{linenomath}
\begin{equation}
 \int_{-\infty}^{+\infty}f_X(x')\,dx' = 
 \int_{V_{X}}f_X(x')\,dx' = 1 \ .
\label{eq:norpdf}
\end{equation}
\end{linenomath}
Two pdf's correspond to the same cdf precisely if they differ only on a
set of Lebesgue measure zero. On the other hand, a cdf of a continuous random
variable is differentiable almost
everywhere on $\mathbb{R}$ (\cite{ste}, \S\,3.2, Theorem\;3.11,
  pp.\,130-131) such that the derivative can be used as a pdf.
\begin{define}
\label{def:scpdf}
Throughout the present discussion,
\begin{linenomath}
\begin{equation}
 f_X(x)\equiv\frac{d}{dx}F_X(x) 
\label{eq:cdf1}
\end{equation}
\end{linenomath}
is assumed.
\end{define}
\begin{define}[Probability distribution]
\label{def:pr}
A function $Pr_{X}:{\mathcal B}\longrightarrow [0,1]$ called
     {\em pro\-ba\-bi\-li\-ty distribution} is defined as
     the image measure of $P$ by the random variable $X$, 
     $Pr_{X}\equiv P\circ X^{-1}$,
     such that $Pr_X(S)=P[X^{-1}(S)]$, where 
     $X^{-1}(S)\in\Sigma$ is the inverse image of a Borel set $S$
     under $X$. A probability distribution over a continuous random 
     variable $X$ is called a {\em continuous probability distribution}.
\end{define}
From the properties of the underlying probability spaces
it follows immediately that 
probability distributions for random variables also conform to the axioms
(\ref{eq:axiom1}-\ref{eq:axiom3}) of probability. Therefore, a scalar 
random variate $X$ on a probability space $(\Omega,\Sigma,P)$ 
generates another probability space 
$(\mathbb{R},{\mathcal B},Pr_X)$ with the Borel algebra 
${\mathcal B}\equiv{\mathcal B}^1$ as underlying $\sigma$-algebra.

Let $X$ and $Y$ be continuous random variables defined on  
$(\Omega,\Sigma,P)$, let there exist 
a function $s$ on $V_X$ such that   $Y= s\circ X$ and 
$y=s(x)$, and let the function $s$ be
differentiable with non-vanishing derivative $s'(x)$ 
on the entire support $V_{X}$ of $X$, such that 
$[s^{-1}(y)]'=[s'(x)]^{-1}$ exists for all $y=s(x)$ with $x\in V_X$.
Then, due to the common probability space $(\Omega,\Sigma,P)$
underlying the spaces $(\mathbb{R},{\mathcal B},Pr_X)$ and
$(\mathbb{R},{\mathcal B},Pr_Y)$,
\begin{linenomath}
$$
\beginpicture
 \setcoordinatesystem units <1 mm,1 mm>
 \put{$(\Omega,\Sigma,P)$} at 0 0
 \put{$(\mathbb{R},{\mathcal B},Pr_X)$} at -19 17
 \put{$(\mathbb{R},{\mathcal B},Pr_Y)$} at 19 17
 \put{$X$} at -7.5 10.5
 \put{$Y$} at 8.5 10.5
 \put{$s$} at 0 19
 \put{$,$} at 35 8
 \arrow <1.5 mm> [.3,.6] from -8 17 to 8 17
 \arrow <1.5 mm> [.3,.6] from -3 4 to -17 14
 \arrow <1.5 mm> [.3,.6] from 3 4 to 17 14
\endpicture
$$
\end{linenomath}
for all $y$ for which $s^{-1}(y)\in V_X$
the cdf $F_Y$ for $Y$ can be expressed in terms of $F_X$ as
\begin{linenomath} 
\begin{equation}
 F_Y(y)=
\begin{cases}
 \hskip 4mm F_X(s^{-1}(y)) \hskip 4mm 
  ; \ [s^{-1}(y)]' > 0 \\
  1-F_X(s^{-1}(y)) \  ; \ [s^{-1}(y)]' < 0
\end{cases} 
,
\label{eq:leminvar}
\end{equation} 
\end{linenomath}
and the pdf for $Y$ is related to the pdf for $X$ as
\begin{linenomath} 
\begin{equation} 
 f_Y(y)=\frac{d}{dy}F_Y(y)=
 f_X(s^{-1}(y))\,\mbox{\large $|$}[s^{-1}(y)]'\mbox{\large $|$} .
 \label{eq:vartr0}
\end{equation} 
\end{linenomath}
The image of $V_X$ under $s$ is contained in $V_Y$, $s(V_X)\subseteq
V_Y$, and the probability distribution $Pr_Y[V_Y-s(V_X)]$
for the relative complement of $V_Y$ and $s(V_X)$ is zero.

The foregoing discussion about the probability distributions
associated to scalar random variables is extended to
multivariate random variables as follows.
\begin{define}[Random vectors]
\label{def:ranvec}
Given a probability space $(\Omega,\Sigma,P)$, a vector
 function $\mathbf{X}=( X_1,\ldots , X_n)$ is called a
 {\em multivariate random variable} (or {\em random vector})
 if $A_{\bi{X}\le\bi{x}}=\{\omega\in \Omega : X_1(\omega)\le
x_1,\ldots,X_n(\omega)\le x_n\}\in\Sigma$, $\forall\,\bi{x}
=(x_1,\ldots ,x_n)\in\mathbb{R}^{n}$.
Every random vector gives rise to a cdf $F_{\bi{X}}(x_1,\ldots
,x_n)$ on the state space $\mathbb{R}^n$ to $[0,1]$ such
that $ F_{\bi{X}}(x_1,\ldots,x_n)=P(A_{\bi{X}\le\bi{x}})$, 
and to a {\em joint} probability distribution 
$Pr_{\bi{X}}\left( S\right)$ 
on the Borel algebra ${\mathcal B}^n$ 
to $[0,1]$, $Pr_{\bi{X}}\left( S\right)\equiv 
P[\bi{X}^{-1}(S)]$, $S\in{\mathcal B}^n$.
Also, as for the scalar random variates, 
a random vector $\bi{X}$ is called continuous if its cdf
can be written as an integral of a pdf $f_{\bi{X}}(x_1,\ldots,x_n)$,
\begin{linenomath} 
\begin{equation*}
\begin{split}
 F_{\bi{X}}(x_1,\ldots,x_n)&=\int_{U_{\bi{X}\le\bi{x}}}
 f_{\bi{X}}(t_1,\ldots,t_n)\,dt_1\dotsm dt_n \\
 &=
 \int_{-\infty}^{x_1}dt_1\dotsm\int_{-\infty}^{x_n}dt_n\, 
 f_{\bi{X}}(t_1,\ldots,t_n) ,
\end{split}
\end{equation*} 
\end{linenomath}
where $U_{\bi{X}\le\bi{x}}\equiv \times_{i=1}^n (-\infty,x_i]$ 
is an infinite $n$-dimensional rectangle in the 
state space $\mathbb{R}^n$, while the transition from a 
$n$-dimensional integral to $n$ iterated
integrals is justified by Fubini's Theorem (see, for
  example, \cite{bar}, Chapter\;10, pp.\,119-120).
\end{define}
Every (joint) probability distribution for a continuous $n$-vector $\bi{X}$ 
can be expressed as an integral
\begin{linenomath} 
\begin{equation*}
 Pr_{\bi{X}}(S)\equiv\int_{S}f_{\bi{X}}(\bi{x})\,d^n\bi{x} 
\ ; \ \ \ \ \forall\,S\in{\mathcal B}^n \ . 
\end{equation*} 
\end{linenomath}

Let $\bi{X}$ and $\bi{Y}$ be $n$-dimensional continuous random variables on a
probability space $(\Omega,\Sigma,P)$, let
$f_{\bi{X}}(\bi{x})$ be a pdf for $\bi{X}$, and let $\bi{s}$ be a
differentiable function on $V_{\bi{X}}$ with non-vanishing Jacobian $\mbox{\large $|$}
\partial_{\bi{x}}\bi{s}(\bi{x})\mbox{\large $|$}$ such that
$\bi{Y}=\bi{s}\circ\bi{X}$.
Then, for all $\bi{y}$ from the image of $V_{\bi{X}}$ under
  $\bi{s}$, the pdf for $\bi{Y}$ reads:
\begin{linenomath} 
\begin{equation}
 f_{\bi{Y}}(\bi{y})\equiv 
 f_{\bi{X}}(\bi{s}^{-1}(\bi{y}))\,
 |\partial_{\bi{y}}\bi{s}^{-1}(\bi{y})| \ .
\label{eq:vartr1}
\end{equation} 
\end{linenomath}
\begin{define}[Marginal distributions]
\label{def:marginal}
Let a random vector $\bi{X}$ be partitioned into 
a random $n$-vector
$\bi{Y}$ and a random $m$-vector $\bi{Z}$, $\bi{X}=\left(\bi{Y},
\bi{Z}\right)$. Then $F_{\mathbf{X}}^{\bi{Y}}(\bi{y}) \equiv 
 F_{\mathbf{X}}(\bi{y},z_1=\infty,\ldots,z_m=\infty)$ and
$ F_{\mathbf{X}}^{\bi{Z}}(\bi{z}) \equiv 
 F_{\mathbf{X}}(y_1=\infty,\ldots,y_n=\infty,\bi{z})$
are called the {\em marginal cdf's} for the
components $\bi{Y}$ and $\bi{Z}$ of the partition
$\left(\bi{Y},\bi{Z}\right)$ of  $\mathbf{X}$, respectively.
Also, pdf's
\begin{linenomath} 
\begin{equation*}
 f^{\bi{Y}}_{\mathbf{X}}(\bi{y}) \equiv 
 \int_{\mathbb{R}^{m}}f_{\mathbf{X}}(\bi{y},\bi{z})\,d^{m}\bi{z}
\end{equation*} 
\end{linenomath}
and
\begin{linenomath} 
\begin{equation*}
 f^{\bi{Z}}_{\mathbf{X}}(\bi{z}) \equiv 
 \int_{\mathbb{R}^{n}}f_{\mathbf{X}}(\bi{y},\bi{z})\,d^{n}\bi{y}
\end{equation*} 
\end{linenomath}
are called the {\em marginal pdf's} for the components $\bi{Y}$
and $\bi{Z}$ of a partition of a continuous random vector 
$\mathbf{X}$, while the corresponding {\em marginal probability
  distributions} are denoted by $Pr_{\bi{X}}^{\bi{Y}}(U)$ and
  $Pr_{\bi{X}}^{\bi{Z}}(S)$, $U\in{\mathcal B}^n$ and $S\in{\mathcal B}^m$. 
\end{define}
Usually, abbreviated notations may be used, e.g., 
$F_{\mathbf{X}}(\bi{y})\equiv F^{\bi{Y}}_{\mathbf{X}}(\bi{y})$ and
$f_{\mathbf{X}}(\bi{z})\equiv f^{\bi{Z}}_{\mathbf{X}}(\bi{z})$.
Since, however, in $F_{\mathbf{X}}(\bi{y})$ 
and in $f_{\mathbf{X}}(\bi{y})$ the
arguments of the functions denote also the functions themselves, 
it should be noted that $F_{\mathbf{X}}(\bi{y})$ and
$f_{\mathbf{X}}(\bi{y})$ are not necessarily the same
functions as $F_{\mathbf{X}}(\bi{z})$ and
$f_{\mathbf{X}}(\bi{z})$, respectively. 
\begin{define}[Conditional probability distributions]
\label{def:condprdis}
Let $(\Omega,\Sigma,P)$ be a probability space and
$\bi{X}=(\bi{Y},\bi{Z}):
\Omega\longrightarrow\mathbb{R}^{n}\times\mathbb{R}^{m}$ a 
$\Sigma$-measurable function that gives rise to a probability
distribution $Pr_{\bi{X}}: {\mathcal B}^{n}\times {\mathcal B}^{m}
\longrightarrow [0,1]$,
let $(\mathbb{R}^n,{\mathcal B}^n,Pr_{\bi{X}}^{\bi{Y}})$ and 
$(\mathbb{R}^m,{\mathcal B}^m,Pr_{\bi{X}}^{\bi{Z}})$ be 
the spaces of the marginal probability distributions for the
components $\bi{Y}$ and  $\bi{Z}$ of the partition $(\bi{Y},\bi{Z})$ 
of $\bi{X}$, and
let  $\mathbf{1}_{\bi{Y}^{-1}(U)}$, $U\in{\mathcal B}^n$, be the 
{\em indicator function} on $\Omega$:
$\mathbf{1}_{\mathbf{Y}^{-1}(U)}(\omega)=1$ for $\omega\in
 \bi{Y}^{-1}(U)$ and $0$ otherwise.
 Then, a function
$\nu_{\bi{1}_{\bi{Y}^{-1}(U)}}:\Sigma'\longrightarrow\mathbb{R}$, $\Sigma'\equiv
\bi{Z}^{-1}({\mathcal B}^m)\subset\Sigma$,
\begin{linenomath} 
\begin{equation*}
 \nu_{\bi{1}_{\bi{Y}^{-1}(U)}}[\bi{Z}^{-1}(S)]\equiv
 \int_{\bi{Z}^{-1}(S)}\!\!\mathbf{1}_{\mathbf{Y}^{-1}(U)}(\omega)\,
 dP(\omega) \ ,
\end{equation*} 
\end{linenomath}
$S\in{\mathcal B}^m$, is a finite measure on $\Sigma'$, and so is finite
the image measure  $\widetilde{\nu}_{\bi{1}_{\bi{Y}^{-1}(U)}}$ of the measure
$\nu_{\bi{1}_{\bi{Y}^{-1}(U)}}$ by $\bi{Z}$, 
$\widetilde{\nu}_{\bi{1}_{\bi{Y}^{-1}(U)}}: {\mathcal B}^m
\longrightarrow \mathbb{R}$, $\widetilde{\nu}_{\bi{1}_{\bi{Y}^{-1}(U)}}
\equiv \nu_{\bi{1}_{\bi{Y}^{-1}(U)}}\circ \bi{Z}^{-1}$. The function
$Pr_{\bi{X}}^{\bi{Y}|\bi{Z}=\bi{z}}(U|\bi{z}):
\mathbb{R}^m\longrightarrow\mathbb{R}$  
called  {\em conditional probability distribution for $\bi{Y}$
given the value $\bi{Z}=\bi{z}$}, is then defined by the set
of functional equations:
\begin{linenomath} 
\begin{equation}
 \widetilde{\nu}_{\bi{1}_{\bi{Y}^{-1}(U)}}(S)=
 \int_{S}Pr_{\bi{X}}^{\bi{Y}|\bi{Z}=\bi{z}}(U|\bi{z})(\bi{z})\,
 dPr_{\bi{X}}^{\bi{Z}}(\bi{z}) \ ,
\label{eq:defcondprdis}
\end{equation} 
\end{linenomath}
while the corresponding {\em conditional cdf} is denoted by
$F_{\bi{X}}^{\bi{Y}|\bi{Z}=\bi{z}}(\bi{y}|\bi{z})$.
\end{define}
The definition of $Pr_{\bi{X}}^{\bi{Y}|\bi{Z}=\bi{z}}(U|\bi{z})$ can 
be interpreted to say that the diagram 
\begin{linenomath}
$$
\beginpicture
 \setcoordinatesystem units <1 mm,1 mm>
 \put{$\Omega$} at 0 0
 \put{$\mathbb{R}$} at 0 18
 \put{$\mathbb{R}^m$} at 22 0 
 \put{$\bi{1}_{\bi{Y}^{-1}(U)}$} at -7 9
 \put{$\bi{Z}^{-1}$} at 12 -3
 \put{$Pr_{\bi{X}}^{\bi{Y}|\bi{Z}=\bi{z}}(U|\bi{z})$} at 25.5 9
 \arrow <1.5 mm> [.3,.6] from 0 3 to 0 15
 \arrow <1.5 mm> [.3,.6] from 18 0  to 3 0
 \arrow <1.5 mm> [.3,.6] from 19 2 to 2 16
\endpicture
$$
\end{linenomath}
is commutative in the average with respect to $Pr_{\bi{X}}^{\bi{Z}}$.
\begin{define}[Conditional pdf]
\label{def:condpdf}
Let $Pr_{\bi{X}}^{\bi{Y}|\bi{Z}=\bi{z}}(U|\bi{z})$ be a solution of  
\eqref{eq:defcondprdis}.
For continuous $\bi{X}$, the system of equations
\begin{linenomath} 
\begin{equation}
 Pr_{\bi{X}}^{\bi{Y}|\bi{Z}=\bi{z}}(U|\bi{z})=\int_{U}\,
 f_{\bi{X}}^{\bi{Y}|\bi{Z}=\bi{z}}(\bi{y}|\bi{z})\,d^n\bi{y}
\label{eq:defcondpdf}
\end{equation} 
\end{linenomath}
for all $ U\in{\mathcal B}^n$, is the defining condition for the 
{\em conditional pdf 
$f_{\bi{X}}^{\bi{Y}|\bi{Z}=\bi{z}}$ for $\bi{Y}$ given $\bi{Z}=\bi{z}$}.
\end{define}
For conditional cdf's and pdf's, abbreviated notations 
$F_{\mathbf{X}}(\bi{y}|\bi{z})\equiv
F_{\mathbf{X}}^{\bi{Y}|\bi{Z}=\bi{z}}(\bi{y}|\bi{z})$ and 
$f_{\mathbf{X}}(\bi{y}|\bi{z})\equiv
f_{\mathbf{X}}^{\bi{Y}|\bi{Z}=\bi{z}}(\bi{y}|\bi{z})$
may again be used.
\begin{theo}
\label{theo:condpdf}
Let $f_{\bi{X}}(\bi{y},\bi{z})$ be a joint pdf for a $(n+m)$-dimensional 
random vector $\bi{X}=(\bi{Y},\bi{Z})$ and let $f_{\bi{X}}(\bi{z})$ be
the marginal pdf for $\bi{Z}$, supported on $V_{\bi{Z}}$. Then, 
\begin{linenomath} 
\begin{equation}
  f_{\bi{X}}(\bi{y}|\bi{z})=\frac{f_{\bi{X}}(\bi{y},\bi{z})}
 {f_{\bi{X}}(\bi{z})}
\label{eq:condpdf}
\end{equation} 
\end{linenomath}
holds true uniquely on $(\mathbb{R}^n-U_0)\times(V_{\bi{Z}}-S_0)$,
where $Pr_{\bi{X}}^{\bi{Z}}(S_0)=\nu_L(U_0)=0$,
$\nu_L(U_0)\equiv\int_{U_0}d^n\bi{y}$. It is said that
$f_{\bi{X}}(\bi{y}|\bi{z})$ is determined uniquely {\em
  $Pr_{\bi{X}}^{\bi{Z}}$-almost everywhere on $V_{\bi{Z}}$} and
{\em $\nu_L$-almost everywhere on $\mathbb{R}^n$}.
\end{theo}
\begin{remark}
\label{rem:conddistr}
{\rm First, the reason for adopting an indirect definition of the
conditional pdf's is that the more direct formulations like, for example,
the approach that is based on the L'H\^{o}pital rule 
(see, for example, \cite{rao}, \S\,1.4,
  pp.\,13-14) and the axiomatization of \cite{ren}, 
do not lead to uniquely defined conditional
pdf's. For a discussion on the resulting inconsistencies see
\cite{rao}, Chapters 3 and 4, pp.\,63-121. Second, below, 
existence of a joint pdf $f_{\bi{X}}(\bi{y},\bi{z})$
is not a necessary condition for existence of the corresponding
conditional pdf's  $f_{\bi{X}}(\bi{y}|\bi{z})$ and $f_{\bi{X}}(\bi{z}|\bi{y})$.
}
\end{remark}

Let there exist a conditional pdf
$f_{\bi{X}}\left(\bi{y},\bi{z}|\bi{t}\right)$, 
$\bi{X}=\left(\bi{Y},\bi{Z},\bi{T}\right)$, and let the
marginal distribution 
\begin{linenomath} 
\begin{equation*}
 f_{\bi{X}}\left(\bi{z}|\bi{t}\right)\equiv 
 \int_{\mathbb{R}^{n_1}}
 f_{\bi{X}}\left(\bi{y},\bi{z}|\bi{t}\right)\,d^{n_1}\bi{y} 
\end{equation*} 
\end{linenomath}
be positive. Then, 
by an iterative application of Definition\;\ref{def:condpdf},
\begin{linenomath} 
\begin{equation}
 f_{\mathbf{X}}(\bi{y}|\bi{z},\bi{t})=
 \frac{f_{\mathbf{X}}(\bi{y},\bi{z}|\bi{t})}
      {f_{\bi{X}}(\bi{z}|\bi{t})}
\ .
\label{eq:condrat1.a}
\end{equation} 
\end{linenomath}
The results of the following example are obtained by sequential 
applications of the {\em product rule} \eqref{eq:condrat1.a}.
\begin{exam}
\label{ex:bayes}
{\rm Let $\bi{X}$ be partitioned into
  $(\bi{Y},\bi{Z},\bi{T},\bi{W})$ and let
there exist conditional pdf's 
$f_{\bi{X}}\left(\bi{y},\bi{t}|\bi{z},\bi{w}\right)$ and
and $f_{\bi{X}}\left(\bi{y},\bi{w}|\bi{z},\bi{t}\right)$.
Then, in an analogy with \eqref{eq:condrat1.a}, for 
$f_{\bi{X}}\left(\bi{t}|\bi{z},\bi{w}\right),  
f_{\bi{X}}\left(\bi{w}|\bi{z},\bi{t}\right) > 0$ there exists a
conditional pdf 
$f_{\bi{X}}\left(\bi{y}|\bi{z},\bi{t},\bi{w}\right)$ such that 
\begin{linenomath} 
\begin{equation*}
 f_{\bi{X}}\left(\bi{y}|\bi{z},\bi{t},\bi{w}\right)=
 \frac{f_{\bi{X}}\left(\bi{y},\bi{t}|\bi{z},\bi{w}\right)}
      {f_{\bi{X}}\left(\bi{t}|\bi{z},\bi{w}\right)} 
=
 \frac{f_{\bi{X}}\left(\bi{y},\bi{w}|\bi{z},\bi{t}\right)}
      {f_{\bi{X}}\left(\bi{w}|\bi{z},\bi{t}\right)} \ .
\label{eq:condrat1}
\end{equation*} 
\end{linenomath}
When, in addition, the marginal pdf's 
$f_{\bi{X}}\left(\bi{y}|\bi{z},\bi{w}\right)$ and
$f_{\bi{X}}\left(\bi{y}|\bi{z},\bi{t}\right)$ are also non-vanishing,
the joint pdf's 
$f_{\bi{X}}\left(\bi{y},\bi{t}|\bi{z},\bi{w}\right)$ and
$f_{\bi{X}}\left(\bi{y},\bi{w}|\bi{z},\bi{t}\right)$
can be further decomposed as
$f_{\bi{X}}\left(\bi{y},\bi{t}|\bi{z},\bi{w}\right)=
 f_{\bi{X}}\left(\bi{y}|\bi{z},\bi{w}\right)\,
 f_{\bi{X}}\left(\bi{t}|\bi{y},\bi{z},\bi{w}\right)$ and 
$ f_{\bi{X}}\left(\bi{y},\bi{w}|\bi{z},\bi{t}\right)=
f_{\bi{X}}\left(\bi{y}|\bi{z},\bi{t}\right)\,
 f_{\bi{X}}\left(\bi{w}|\bi{y},\bi{z},\bi{t}\right)$,
such that
\begin{linenomath} 
\begin{equation}
\begin{split}
 f_{\bi{X}}\left(\bi{y}|\bi{z},\bi{t},\bi{w}\right)&=
 \frac{f_{\bi{X}}\left(\bi{y}|\bi{z},\bi{w}\right)\,
 f_{\bi{X}}\left(\bi{t}|\bi{y},\bi{z},\bi{w}\right)}
      {f_{\bi{X}}\left(\bi{t}|\bi{z},\bi{w}\right)} \\
&=
 \frac{f_{\bi{X}}\left(\bi{y}|\bi{z},\bi{t}\right)\,
       f_{\bi{X}}\left(\bi{w}|\bi{y},\bi{z},\bi{t}\right)}
      {f_{\bi{X}}\left(\bi{w}|\bi{z},\bi{t}\right)} \ .
\end{split}
\label{eq:bayes0}
\end{equation} 
\end{linenomath}
In the same way,
\begin{linenomath} 
\begin{equation}
\begin{split}
 f_{\bi{X}}\left(\bi{y}|\bi{t},\bi{w}\right)
&=
 \frac{f_{\bi{X}}\left(\bi{y}|\bi{w}\right)\,
       f_{\bi{X}}\left(\bi{t}|\bi{y},\bi{w}\right)}
      {f_{\bi{X}}\left(\bi{t}|\bi{w}\right)} \\
&=
 \frac{f_{\bi{X}}\left(\bi{y}|\bi{t}\right)\,
       f_{\bi{X}}\left(\bi{w}|\bi{y},\bi{t}\right)}
      {f_{\bi{X}}\left(\bi{w}|\bi{t}\right)} 
\end{split}
\label{eq:bayes01}
\end{equation} 
\end{linenomath}
is obtained when $\bi{X}$ is partitioned into $(\bi{Y},\bi{T},\bi{W})$.
}
\end{exam}
\begin{exam}[Transformations of conditional pdf's]
\label{ex:transcond}
{\rm Let $\mathbf{X}=(\bi{X}_1,\bi{X}_2)$ be a continuous
  ($n_1+n_2$)-dimensional random variable and 
 $f_{\mathbf{X}}(\bi{x}_1|\bi{x}_2)$ be a conditional pdf.
 Let, in addition, $\bi{s}:
 V_{\bi{X}}
\longrightarrow \mathbb{R}^{n_1}\times\mathbb{R}^{n_2}$ 
be a differentiable function function such that 
 $\bi{Y}\equiv(\bi{Y}_1,\bi{Y}_2)=\bi{s}\circ\bi{X}=
 (\bi{s}_1\circ \bi{X}_1,\bi{s}_2\circ \bi{X}_2)$ and
 that the Jacobian $\mbox{\large $|$}\partial_{\bi{x}}\bi{s}(\bi{x})\mbox{\large $|$}
 = \mbox{\large $|$}\partial_{\bi{x}_1}\bi{s}_1(\bi{x}_1)\mbox{\large $|$}\,
\mbox{\large $|$}\partial_{\bi{x}_2}\bi{s}_2(\bi{x}_2)\mbox{\large $|$}$ does not 
 vanish on the entire support $V_{\bi{X}}$ of $\bi{X}$. 
 For $f_{\mathbf{X}}(\bi{s}_2^{-1}(\bi{y}_2))>0$, 
 equations \eqref{eq:vartr0} and \eqref{eq:vartr1}
 applied to the conditional pdf $f_{\mathbf{Y}}(\bi{y}_1|\bi{y}_2)
 = f_{\mathbf{Y}}(\bi{y}_1,\bi{y}_2)/f_{\mathbf{Y}}(\bi{y}_2)$ 
 then yield
\begin{linenomath} 
\begin{equation}
\begin{split}
 f_{\mathbf{Y}}(\bi{y}_1|\bi{y}_2)
      &=
 \frac{f_{\mathbf{X}}(\bi{s}_1^{-1}(\bi{y}_1),\bi{s}_2^{-1}(\bi{y}_2))}
      {f_{\mathbf{X}}(\bi{s}_2^{-1}(\bi{y}_2))}\,
 \mbox{\large $|$}\partial_{\bi{y}_1}\bi{s}_1^{-1}(\bi{y}_1)\mbox{\large $|$} 
\\ 
&=
 f_{\mathbf{X}}(\bi{s}_1^{-1}(\bi{y}_1)|\bi{s}_2^{-1}(\bi{y}_2))
\mbox{\large $|$}\partial_{\bi{y}_1}\bi{s}_1^{-1}(\bi{y}_1)\mbox{\large $|$}
\ .
\end{split}
\label{eq:condtrans0}
\end{equation} 
\end{linenomath}
}
\end{exam}

During the present discussion we allow for a possibility that a 
conditional pdf $f_{\mathbf{X}}(\bi{x}_1|\bi{x}_2)$ exists even when the 
corresponding joint pdf $f_{\mathbf{X}}(\bi{x}_1,\bi{x}_2)$ does not exist.
When $f_{\mathbf{X}}(\bi{x}_1,\bi{x}_2)$ does not exist,
however, the transformation \eqref{eq:condtrans0} of the
conditional pdf that is induced by the transformation of the random
vector, ceased to be uniquely determined. In order to dismiss this
ambiguity, the following definition, motivated by the 
preceding example, is adopted.
\begin{define}[Transformations of conditional pdf's]
\label{def:transcond}
Let there exist a conditional pdf $f_{\mathbf{X}}(\bi{x}_1|\bi{x}_2)$,
$\bi{X}=(\bi{X}_1,\bi{X}_2)$ and $\bi{x}=(\bi{x}_1,\bi{x}_2)$, and let 
a function $\bi{s}: (\bi{x}_1,\bi{x}_2)\longrightarrow(\bi{s}_1(\bi{x}_1), 
 \bi{s}_2(\bi{x}_2))\equiv(\bi{y}_1,\bi{y}_2)$
     be one-to-one and with non-vanishing Jacobian 
     $|\partial_{\bi{x}_1}\bi{s}_1(\bi{x}_1)|$ 
     on the entire support $V_{\bi{X}_1|\bi{x}_2}$ of
     $f_{\mathbf{X}}(\bi{x}_1|\bi{x}_2)$.
     Then, the conditional pdf
     $f_{\mathbf{Y}}(\bi{y}_1|\bi{y}_2)$,
     $\bi{Y}\equiv(\bi{s}_1\circ\bi{X}_1, \bi{s}_2\circ
     \bi{X}_2)\equiv(\bi{Y}_1,\bi{Y}_2)$,
     is defined as
\begin{linenomath} 
\begin{equation}
 f_{\mathbf{Y}}(\bi{y}_1|\bi{y}_2)\equiv
 f_{\mathbf{X}}(\bi{s}_1^{-1}(\bi{y}_1)|\bi{s}_2^{-1}(\bi{y}_2))
\mbox{\large $|$}\partial_{\bi{y}_1}\bi{s}_1^{-1}(\bi{y}_1)\mbox{\large $|$}
\ ,
\label{eq:condtrans1}
\end{equation} 
\end{linenomath}
where $\bi{s}_{1,2}^{-1}$ are the inverse functions of $\bi{s}_{1,2}$.
\end{define}
%
%

\subsection{Parametric families of probability distributions}
\label{ss:parfam}

The term {\em parametric family} is used to describe a collection 
$I=\{Pr_{I,\bs{\theta}}: \bs{\theta}\in V_{\mathbf{\Theta}}\}$
of probability distributions that differ only in the value of a
(possibly multi-di\-men\-sio\-nal) {\em parameter}, say $\mathbf{\Theta}$, i.e.,
a value $\bs{\theta}$ of $\mathbf{\Theta}$ determines a unique distribution
within $I$. Therefore, a probability distribution 
for a random $n$-vector $\bi{X}$,
$Pr_{\bi{X}}\left(S_{\bi{X}}\right)$,
$S_{\bi{X}}\in{\mathcal B}^n$, that belongs to a particular parametric
family $I$, is denoted by
$Pr_{I,\bs{\theta}}\left(S_{\bi{X}}\right)$, 
whereas $F_{I,\bs{\theta}}\left(\bi{x}\right)$ stands for the
corresponding cdf. Likewise, $f_{I,\bs{\theta}}\left(\bi{x}\right)$
denotes a unique pdf within a parametric family $I$ of continuous
probability distributions. A continuous probability distribution from 
a parametric family $I$ is supported on a set
$V_{\bi{X}}=V_{\bi{X}}(\bs{\theta})$ that may, in general, depend on
the value $\bs{\theta}$ of the parameter, while the range
$V_{\mathbf{\Theta}}\subseteq\mathbb{R}^{m}$ of admissible values of
$\mathbf{\Theta}$ is called a {\em parameter space}. In the present
article, every considered parametric family is assumed to be {\em
  identifiable}: $Pr_{I,\bs{\theta}_1}\neq Pr_{I,\bs{\theta}_2}$ for
$\bs{\theta}_1\neq\bs{\theta}_2$, $\bs{\theta}_{1,2}\in V_{\bi{\Theta}}$.

\begin{exam}[Reparameterization]
\label{ex:repar}
{\rm Let $f_{I,\bs{\theta}}\left(\bi{x}\right)$ be a pdf for a random
$n$-vector $\bi{X}$ from a parametric family $I$ and let $\bi{s}$ be a
one-to-one Borel function onto ${\mathbb R}^n$ such that the Jacobian
$\mbox{\large $|$}\partial_{\bi{x}}\bi{s}(\bi{x})\mbox{\large $|$}$ does not vanish 
anywhere on the support $V_{\bi{X}}(\bs{\theta})$ of $\bi{X}$. Then, 
according to \eqref{eq:vartr1},
$ f_{I',\bs{\theta}}(\bi{y})=
 f_{I,\bs{\theta}}(\bi{s}^{-1}(\bi{y}))\,
 \mbox{\large $|$}\partial_{\bi{y}}\bi{s}^{-1}(\bi{y})\mbox{\large $|$}$,
where $\bi{Y}\equiv\bi{s}\circ\bi{X}$, $\bi{y}\equiv\bi{s}(\bi{x})$
and $\bi{s}^{-1}$ is the inverse function of $\bi{s}$, while
indices $I$ and $I'$ indicate that probability distributions for
$\bi{X}$ and $\bi{Y}$ in general belong to different (but isomorphic)
parametric families. 
Let, in addition, $\bar{\bi{s}}$ be a one-to-one function on
the parameter space $V_{\mathbf{\Theta}}\subseteq\mathbb{R}^m$, such
that  $\bs{\lambda}\equiv\bar{\bi{s}}(\bs{\theta})$. Then,
$f_{I',\bs{\theta}}\left(\bi{y}\right)$ can be
reparameterized as
\begin{linenomath} 
\begin{equation}
\begin{split}
 f_{I'',\bs{\lambda}}\left(\bi{y}\right)&=
 f_{I',\bar{\bi{s}}^{-1}(\bs{\lambda})}\left(\bi{y}\right)\\
&=
 f_{I,\bar{\bi{s}}^{-1}(\bs{\lambda})}(\bi{s}^{-1}(\bi{y}))\,
 \mbox{\large $|$}\partial_{\bi{y}}\bi{s}^{-1}(\bi{y})\mbox{\large $|$}\ ,
\end{split}
\label{eq:repar}
\end{equation} 
\end{linenomath}
where $\bar{\bi{s}}^{-1}$ is the inverse function of $\bar{\bi{s}}$.
}
\end{exam}
There is a complete analogy between the transformation
\eqref{eq:repar} and the transformations \eqref{eq:condtrans0} and
\eqref{eq:condtrans1}, such that every probability distribution from 
a parametric family can be regarded as a conditional distribution, 
i.e., as a distribution that is conditional upon the value of the
parameter. Accordingly, we define
$ F_{I}\left(\bi{x}|\bs{\theta}\right)\equiv 
 F_{I,\bs{\theta}}\left(\bi{x}\right)$ and
$Pr_{I}\left(S_{\mathbf{\Theta}=\bs{\theta},\bi{X}}|\bs{\theta}\right)\equiv 
  Pr_{I,\bs{\theta}}\left(S_{\bi{X}}\right)$, 
$S_{\mathbf{\Theta}=\bs{\theta},\bi{X}}\in{\mathcal B}^n$,
and, for continuous $\bi{X}$,
\begin{linenomath} 
\begin{equation}
 f_{I}\left(\bi{x}|\bs{\theta}\right)\equiv 
 f_{I,\bs{\theta}}\left(\bi{x}\right) 
\label{eq:directpdf}
\end{equation} 
\end{linenomath}
for all $\bi{x}\in\mathbb{R}^n$ and $\bs{\theta}\in
V_{\mathbf{\Theta}}\subseteq\mathbb{R}^m$. 

$F_{I}\left(\bi{x}|\bs{\theta}\right)$, 
$Pr_{I}\left(S_{\mathbf{\Theta}=\bs{\theta},\bi{X}}|\bs{\theta}\right)$
and $f_{I}\left(\bi{x}|\bs{\theta}\right)$ are underlain
by a probability space
$(\Omega_{\bs{\theta}},\Sigma_{\bs{\theta}},P)$ and by a 
$(m+n)$-dimensional random
variable $(\bi{\Theta},\bi{X}): \Omega_{\bs{\theta}}\longrightarrow
\left(\bs{\theta},\mathbb{R}^n\right)$ for all $\bs{\theta}\in
V_{\bi{\Theta}}$, where every state space
$\left(\bs{\theta},\mathbb{R}^n\right)$ is a slice on
$V_{\bi{\Theta}}\times\mathbb{R}^n$ that corresponds to a particular
value $\bs{\theta}$ of a $m$-dimensional parameter $\bi{\Theta}$ of the
family $I$. The probability distributions 
$Pr_{I}\left(S_{\mathbf{\Theta}=\bs{\theta},\bi{X}}|\bs{\theta}\right)$ on
Borel $\sigma$-algebras ${\mathcal B}^n$ on such slices are called
{\em direct probability distributions} and represent the first
step towards a unified approach to random variables and parameters from
parametric families. The second step is made in
Section\;\ref{sec:inverse}, where the notion of the inverse
probability distribution is introduced.
\begin{remark}
\label{rem:doublenotation}
{\rm The results of Subsections\;\ref{ss:parfam} and \ref{ss:invariance}
are independent of the preceding definitions. 
The only reason to define $F_{I}\left(\bi{x}|\bs{\theta}\right)$
and $f_{I}\left(\bi{x}|\bs{\theta}\right)$ already at this stage is
to avoid unnecessary duplications in notation.
}
\end{remark}
\begin{define}[Independent random variables]
\label{def:independent}
When 
$f_{I}(\bi{x}|\bi{y},\bs{\theta})= f_I(\bi{x}|\bs{\theta})$ 
and $f_I(\bi{y}|\bi{x},\bs{\theta})= f_I(\bi{y}|\bs{\theta})$, 
the components $\bi{X}$ and $\bi{Y}$ of a continuous random vector 
$(\bi{X},\bi{Y})$ are called {\em independent random variables}.
When, in addition, $f_I(\bi{x}|\bs{\theta})$ and 
$f_I(\bi{y}|\bs{\theta})$ are the same functions, the variables 
$\bi{X}$ and $\bi{Y}$ are 
said to have {\em identical probability distribution}.
\end{define}
When the components $\bi{X}$ and $\bi{Y}$ of a random vector 
$(\bi{X},\bi{Y})$ are independent random variables and the joint pdf
$f_I\left(\bi{x},\bi{y}|\bs{\theta}\right)$ exists, the latter
can be written as 
$f_I\left(\bi{x},\bi{y}|\bs{\theta}\right)=
 f_I(\bi{x}|\bs{\theta})\,
 f_I(\bi{y}|\bs{\theta})$.
\begin{define}[Location and scale parameters]
\label{def:locdis}
Suppose a cdf for a scalar random variable $X$ from a parametric family
$I$ is of the form
\begin{linenomath} 
\begin{equation}
 F_I(x|\mu,\sigma)=
 \Phi\!\left(\frac{x-\mu}{\sigma}\right) \ ,
\label{eq:deflocscale}
\end{equation} 
\end{linenomath}
where $\mu$ is a realization of 
the first component of a two-dimensional parameter 
$\mathbf{\Theta}=(\Theta_1,\Theta_2)$, whereas $\sigma$ is a realization of
its second component. Then, $\Theta_1$ is called a {\em location
parameter} and $\Theta_2$ is called a {\em scale parameter}, while
$V_{\mathbf{\Theta}}=\mathbb{R}\times\mathbb{R}^+$.
\end{define}
When probability distributions from a location-scale family $I$ are
continuous, on the support $V_{X}(\mu,\sigma)$ of a distribution from
the family the appropriate pdf is of the form  
\begin{linenomath} 
\begin{equation}
 f_I(x|\mu,\sigma)=\frac{d}{dx}F_I(x|\mu,\sigma)=
 \frac{1}{\sigma}\,\phi\!\left(\frac{x-\mu}{\sigma}\right) \ ,
\label{eq:defpdflocscale}
\end{equation} 
\end{linenomath}
where $\phi(x)\equiv\Phi'(x)$.
Except for $x=\mu$, every pdf \eqref{eq:defpdflocscale} from 
a location-scale family can be written as a sum 
\begin{linenomath} 
\begin{equation*}
 f_{I}(x|\mu,\sigma)=c_+\,f_{I^+}(x|\mu,\sigma)+c_-\,f_{I^-}(x|\mu,\sigma)
 \ ,
\end{equation*} 
\end{linenomath}
where
\begin{linenomath} 
\begin{equation*}
 c_+\,f_{I^+}(x|\mu,\sigma)\equiv
 \begin{cases}
 \hskip 8mm 
  0
  \hskip 8mm ; \frac{x-\mu}{\sigma} \le 0 \\
 f_{I}(x|\mu,\sigma) \ ; \frac{x-\mu}{\sigma} > 0 \ 
\end{cases}
\end{equation*} 
\end{linenomath}
and
\begin{linenomath} 
\begin{equation*}
 c_-\,f_{I^-}(x|\mu,\sigma)\equiv
 \begin{cases} 
 f_{I}(x|\mu,\sigma) \ ; \frac{x-\mu}{\sigma} < 0 \\
 \hskip 8mm 
  0
  \hskip 8mm ; \frac{x-\mu}{\sigma} \ge 0 
\end{cases}
,
\end{equation*} 
\end{linenomath}
while
\begin{linenomath} 
\begin{equation*}
 c_+\equiv\int_{0}^{\infty}\phi(u)\,du \ \ \ \text{and} \ \ \ 
 c_-\equiv\int_{-\infty}^{0}\phi(u)\,du \ .
\end{equation*} 
\end{linenomath}
For $c_{\pm}>0$, there exist pdf's  
$f_{I^{\pm}}(x|\mu,\sigma)\equiv f_{I}(x|\mu,\sigma)/c_{\pm}$, which
can be further reduced to
\begin{linenomath} 
\begin{equation}
\begin{split}
 f_{{I^{\pm}}'}(y|\lambda_1,\lambda_2=1)&=
 \frac{1}{\lambda_2}\,e^{\frac{y-\lambda_1}{\lambda_2}}\,
 \phi\!\left(\pm e^{\frac{{y-\lambda_1}}{\lambda_2}}\right)\\
&\equiv
 \widetilde{\phi}_{\pm}(y-\lambda_1) \ ,
\end{split}
\label{eq:sc2loc}
\end{equation} 
\end{linenomath}
where $y\equiv\ln{\{\pm(x-\mu)\}} $ and
$\lambda_1\equiv\ln{\sigma}$. That is, every scale parameter for
a location-scale family $I^{\pm}$ is reducible to a location
parameter for a parametric family ${I^{\pm}}'$.
%
\subsection{Invariant families of probability distributions}
\label{ss:invariance}

Let $G=\{a,b,c,\ldots\}$ be a group whose unit element is denoted
by $e$ and let $\bi{l}$ be a function on $G\times\mathbb{R}^{n}$ to
$\mathbb{R}^n$ satisfying $\bi{l}(e,\bi{x})=\bi{x}$,
$\forall\bi{x}\in\mathbb{R}^{n}$ and $\bi{l}(a\circ b,\bi{x})=
\bi{l}[a,\bi{l}(b,\bi{x})]$, 
$\forall a,b\in G$ and
$\forall \bi{x}\in\mathbb{R}^n$. Such a function
specifies $G$ {\em acting on the left of} $\mathbb{R}^n$ and a group
$\mathcal{G}=\{\bi{g}_a:a\in G\}$ of functions $\bi{g}_a : 
\mathbb{R}^n\longrightarrow \mathbb{R}^n$, 
$\bi{g}_a(\bi{x})\equiv\bi{l}(a,\bi{x})$. A composition of 
$\bi{g}_a,\bi{g}_b\in\mathcal{G}$ corresponds to the composition of $a,b\in G$,
$\bi{g}_a[(\bi{g}_b(\bi{x})]\equiv (\bi{g}_a\circ \bi{g}_b)(\bi{x})=
\bi{g}_{a {\scriptscriptstyle{\circ} } b }(\bi{x})$, $\bi{g}_e$ is
the unit element in $\mathcal{G}$ and $\bi{g}_{a^{-1}}=\bi{g}_a^{-1}$, $\forall
a\in G$ (see, for example, \cite{eat}, \S\,2.1, pp.\,19-20).
\begin{define}[Invariant family]
\label{def:invariance}
Let $F_I(\bi{x}|\bs{\theta})$ be a cdf from a parametric
family $I$, let there exist a group $G$ and a function
$\bi{l}: G\times\mathbb{R}^n\longrightarrow\mathbb{R}^n$
specifying both an action of $G$ on the
left of the state space $\mathbb{R}^n$ of the random  
$n$-vector $\bi{X}$ and a group $\mathcal{G}=\{\bi{g}_a:a\in G\}$,
$\bi{g}_a(\,\cdot\,)\equiv\bi{l}(a,\,\cdot\,)$,
and let $\bi{Y}\equiv \bi{g}_a\circ\bi{X}$.
In addition, let for every
$\bi{g}_a\in{\mathcal G}$ and every $\bs{\theta}\in
V_{\mathbf{\Theta}}$ there exist
a transformation $\bar{\bi{g}}_a:\bs{\theta}\longrightarrow 
 \bar{\bi{g}}_a(\bs{\theta})\equiv\bs{\lambda}$, such that
\begin{linenomath} 
\begin{equation}
 F_{I'}(\bi{y}|\bs{\lambda})= F_{I}(\bi{y}|\bs{\lambda}) \ .
\label{eq:definvar}
\end{equation} 
\end{linenomath}
where $\bi{y}\equiv\bi{g}_a(\bi{x})$. The family $I$ is then said to be 
{\em invariant under the group} ${\mathcal G}$ (or ${\mathcal G}${\em
  -invariant} or {\em invariant under the action of the group} $G$) .
\end{define}
Given a ${\mathcal G}$-invariant parametric family $I$, the set 
$\bar{\mathcal G}\equiv \{\bar{\bi{g}}_a:\bi{g}_a\in{\mathcal
G}\}$ of the corresponding transformations on the parameter 
space is also a group (\cite{fer}, \S\,4.1, Lemma\;1, 
pp.\;144-145), usually referred to as the {\em induced group} 
(\cite{ken2}, \S\,23.10, p.\;300).

Let elements of a group $G$ be defined by the values of $n$ continuous
real parameters (or {\em coordinates}), e.g.,
$a=\gamma(a_1,\ldots,a_n)$ with $\gamma$ being a function on (a subset
of) $\mathbb{R}^n$ to $G$. The coordinates are essential in the sense
that the group elements cannot be distinguished by any number of
coordinates smaller
than the {\em dimension} $n$ of the group $G$. Since, by definition, every
group is closed under composition of its elements, $a\circ b=c\in G$,
$\forall a,b\in G$, the coordinates of $c$ are expressible
as functions of the coordinates of $a$ and $b$, 
$c_i=c_i(a_1,\ldots,a_n;b_1,\ldots,b_n)$, $i=1,\ldots,n$. 
\begin{exam}[One-dimensional groups] 
\label{ex:one-dim-gr}
{\rm Coordinates of elements 
of a one-di\-men\-sio\-nal group $G$ also form a group
$\widetilde{G}\subseteq\mathbb{R}$ with $a_1\circ b_1\equiv
c_1(a_1,b_1)$ being the corresponding group operation in 
$\widetilde{G}$. Therefore, since $G$ and $\widetilde{G}$ are isomorphic,
no generality is lost if $a=\gamma(a_1)=a_1$ is assumed.
}
\end{exam}
When coordinates $c_i$ of an element $c=a\circ b$ of a 
$n$-dimensional group $G$ are smooth (i.e., $C^{\infty}$) functions of the
parameters of $a$ and $b$, $G$ is called {\em Lie group}.
\begin{exam}[Invariance of location-scale families]
\label{ex:invlocscale}
{\rm $G= \mathbb{R}\times\mathbb{R}^+$ is a two-di\-men\-sio\-nal Lie group
  for the operations 
\begin{linenomath} 
\begin{equation}
 a\circ b= (a_2b_1+a_1,a_2b_2) \ .
\label{eq:operlocscale}
\end{equation} 
\end{linenomath}
  Every location-scale family $I$ \eqref{eq:deflocscale}
   of continuous probability distributions 
   is invariant under the group 
\begin{linenomath} 
\begin{equation}
 {\mathcal G}=\{g_a: X \longrightarrow a_2X+a_1 \ ;\ 
 (a_1,a_2)\in G\} \ , 
\label{eq:invtransx}
\end{equation} 
\end{linenomath}
with
\begin{linenomath} 
\begin{equation}
\bar{\bi{g}}_{a}: (\Theta_1,\Theta_2) \longrightarrow 
(a_2\Theta_1+a_1,a_2\Theta_2)
\label{eq:invtransms}
\end{equation} 
\end{linenomath}
being the corresponding transformations from the induced group
$\bar{\mathcal G}$. The family $I$ is also invariant under
two one-dimensional subgroups of the group ${\mathcal G}$: under the group
${\mathcal G}_{\times}=\{g_a:X\longrightarrow a X; a\in\mathbb{R}^+\}$, 
$\mathbb{R^+}$ is 
a one-dimensional Lie group for multiplication and 
$\{\bar{\bi{g}}_{a}: (\Theta_1,\Theta_2) \longrightarrow (a\Theta_1,a\Theta_2);
g_a\in{\mathcal G}_{\times}\}$ is the group induced by ${\mathcal G}_{\times}$,
and under
the group ${\mathcal G}_+=\{g_a:X\to a + X; a\in\mathbb{R}\}$,
with $\mathbb{R}$ being a one-dimensional Lie group for summation
and with $\{\bar{\bi{g}}_{a}: (\Theta_1,\Theta_2) \longrightarrow 
(a+\Theta_1,\Theta_2); g_a\in{\mathcal G}_+\}$ being the corresponding
induced group.

Similarly, a family of continuous probability
distributions for random vectors that
consist of two independent scalar random variables $X_1$ and $X_2$, both
belonging to the same location-scale family $I$, is invariant 
under 
${\mathcal G}=\{\bi{g}_a: (X_1,X_2)\longrightarrow
(a_2X_1+a_1,a_2X_2+a_1) \ ;
\ (a_1,a_2)\in \mathbb{R}\times\mathbb{R}^+\}$,
while the corresponding transformations from the induced group are again
\eqref{eq:invtransms}. 
}
\end{exam}
\begin{lemm}
\label{lem:trivial}
Let $G$ be a one-dimensional Lie group, let a function
$l(a,x):G\times\mathbb{R}\longrightarrow\mathbb{R}$ give rise to a
group ${\mathcal G}=\{g_a:a\in G\}$ of transformations 
$g_a:\mathbb{R}\longrightarrow\mathbb{R}$, 
and let $l(a,x)$ be differentiable both in $a$ and in $x$, 
$\forall a\in G$ and $\forall x\in\mathbb{R}$.
Then, for all $x\in\mathbb{R}$ for which
\begin{linenomath} 
\begin{equation}
 \partial_{a}l(a^{-1},x)\mbox{\large $|$}_{a=e} 
\label{eq:prtder1}
\end{equation} 
\end{linenomath}
vanishes, all group transformations are {\em trivial}, i.e.,
$g_a(x) = x$ for all $g_a\in{\mathcal{G}}$.
\end{lemm}
Clearly, if \eqref{eq:prtder1} vanishes for all real $x$, 
then the action of
the group $G$ on the entire real axis is trivial:
$g_a(x) = x$ for every $g_a\in{\mathcal{G}}$ and for all
$x\in\mathbb{R}$.
\begin{lemm}
\label{lem:extheo4}
Suppose a probability distribution for a continuous scalar random variable 
$X$ belongs to a family $I$ of parametric 
distributions that is invariant under the action of a one-dimensional 
Lie group $G$. Let, in addition, the left actions $l(a,x)$ and
$\bar{l}(a,\lambda)$ be differentiable in $a$, $x$ and $\lambda$ 
for all $a\in G$, 
$x\in V_X(\lambda)$ and $\lambda\in V_{\Lambda}(x)$,
let the action of the group $G$ not be identically trivial on the entire 
support $V_{X}(\lambda)$, and let the cdf
for $X$, $F_I(x|\lambda)$, be differentiable in $\lambda$
(differentiability in $x$ is guaranteed by Definition\;\ref{def:scpdf}). 
Then, the partial derivative
\begin{linenomath} 
\begin{equation}
  \partial_{a} \bar{l}(a^{-1},\lambda)\mbox{\large $|$}_{a=e}
\label{eq:prtder2}
\end{equation} 
\end{linenomath}
does not vanish anywhere on the space $V_{\Lambda}$ of the (scalar)
parameter $\Lambda$ of the family $I$.
\end{lemm}
Furthermore, for a continuous scalar random variable $X$ whose
probability distribution belongs to a family $I$ of parametric 
distributions that is invariant under the action of a group $G$,
equation \eqref{eq:leminvar} reduces to
\begin{linenomath} 
\begin{equation}
 F_I\left(x|\lambda\right) =
 \begin{cases}
 \hskip 3.2mm F_I(l(a^{-1},x)|\bar{l}(a^{-1},\lambda))
 \hskip 3.2mm ; [g_a^{-1}(x)]' > 0 \\
 1-   F_I(l(a^{-1},x)|\bar{l}(a^{-1},\lambda)); 
 [g_a^{-1}(x)]' < 0 
\end{cases}
\hskip -4mm,
\label{eq:reslemm1}
\end{equation} 
\end{linenomath}
$a\in G$. On the subspace $\widetilde{V}_X\subseteq V_X(\lambda)$ 
with non-vanishing 
derivatives \eqref{eq:prtder1}, 
derivatives $\partial_a
\bar{l}(a^{-1},\lambda)\mbox{\large $|$}_{a=e}$ are
non-zero by Lemma\;\ref{lem:extheo4}. Then, for $x\in
\widetilde{V}_X$, differentiating 
\eqref{eq:reslemm1} with respect to $a$ and setting afterwards
$a=e$ yields
\begin{linenomath} 
\begin{equation}
\partial_{x}F_I(x|\lambda)\,
 \partial_{\lambda}H(x,\lambda) - 
 \partial_{\lambda}F_I(x|\lambda)\,
 \partial_{x}H(x,\lambda) = 0 \ , 
\label{eq:qpde}
\end{equation} 
\end{linenomath}
where
\begin{linenomath} 
\begin{equation}
 H(x,\lambda)\equiv s(x)-\bar{s}(\lambda) 
\label{eq:defG}
\end{equation} 
\end{linenomath}
and
\begin{linenomath} 
\begin{equation}
[s'(x)]^{-1} \equiv \partial_{a}h(a^{-1},x)\mbox{\large $|$}_{a=e}
\label{eq:defhk1}
\end{equation} 
\end{linenomath}
and
\begin{linenomath} 
\begin{equation}
 [\bar{s}'(\lambda)]^{-1} \equiv 
\partial_{a}\bar{h}(a^{-1},\lambda)\mbox{\large $|$}_{a=e} \ .
\label{eq:defhk2}
\end{equation} 
\end{linenomath}
\begin{lemm}
\label{lem:qpde}
The cdf $F_I(x|\lambda)$ that solves
the functional equation \eqref{eq:qpde}
is a differentiable function of a single variable $H(x,\lambda)$,
\begin{linenomath} 
\begin{equation}
 F_I(x|\lambda) = \Phi[H(x,\lambda)] \ .
\label{eq:lemqpde}
\end{equation} 
\end{linenomath}
\end{lemm}
Consequently, the cdf $F_I(x|\lambda)$ from a parametric family $I$
that is invariant under the action of a one-dimensional Lie 
group can be written as
\begin{linenomath} 
\begin{equation}
 F_I(x|\lambda) = \Phi[s(x)- \bar{s}(\lambda)] =
                 \Phi(y-\mu) \ ,
\label{eq:cdfxth}
\end{equation} 
\end{linenomath}
where $y\equiv s(x)$ and $\mu\equiv \bar{s}(\lambda)$ have been
introduced. Then, by equation \eqref{eq:leminvar}, the cdf for  
the continuous random variable $Y\equiv s\circ X$ 
is of the form
\begin{linenomath} 
\begin{equation*}
 F_{I'}(y|\mu,\sigma=1) = 
 \begin{cases}
\Phi(y-\mu) \ ; \left[s^{-1}(y)\right]' > 0 \\
\widetilde{\Phi}(y-\mu) \ ; \left[s^{-1}(y)\right]' < 0
 \end{cases}
\ ,
\end{equation*} 
\end{linenomath}
where $\widetilde{\Phi}(y-\mu) \equiv 1 - \Phi(y-\mu)$.
That is, the probability distribution for the continuous random variable
$Y$ belongs to a location-scale family $I'$ with $\sigma=1$
(recall equation \eqref{eq:deflocscale}), and the above reasoning can 
be summarized as
\begin{theo}
\label{theo:existence}
Let $X$ be a continuous scalar random variable whose probability
distribution belongs to a ${\mathcal G}$-invariant parametric family
$I$, where ${\mathcal G}=\{g_a:a\in G\}$ is underlain by a
one-dimensional Lie group $G$.
Let, in addition, $g_a(x)$ be differentiable for all
$x\in\mathbb{R}$ and let the 
cdf $F_I(x|\lambda)$ for $X$ be
differentiable in $\lambda$. Then, on the subspace 
$\widetilde{V}_X\subseteq V_X(\lambda)$ with non-vanishing 
derivatives \eqref{eq:prtder1}, $X$ is reducible by a one-to-one
transformation $s$ \eqref{eq:defhk1} to a continuous 
random variable $Y\equiv s\circ X$ whose probability distribution is from a 
location-scale family \eqref{eq:deflocscale} with 
$\sigma=1$ and $\mu\equiv\bar{s}(\lambda)$, where
$\bar{s}$ is defined via \eqref{eq:defhk2}.
\end{theo}
\begin{remark}
\label{rem:theoexistence}
{\rm In the sequel (Proposition\;\ref{theo:existenceprime}) we shall further 
demonstrate that for realizations $x\in V_X(\lambda)-\widetilde{V}_X$ with
vanishing derivative \eqref{eq:prtder1}, a pdf cannot be assigned 
to the inferred parameter of the family $I$.
}
\end{remark}

Let a continuous random variable $\bi{X}$ with a pdf
$f_I(\bi{x}|\bs{\theta})$ 
belong to a parametric family $I$ that is invariant under a group 
${\mathcal G}$ of differentiable 
transformations $\bi{g}_{a}$ with non-vanishing Jacobian 
$\mbox{\large $|$}\partial_{\bi{x}}\bi{g}_{a}(\bi{x})\mbox{\large $|$}$
on the entire support $V_{\bi{X}}(\bs{\theta})$ of the distribution
for $\bi{X}$. Then, equation
\eqref{eq:vartr1} applies which, when combined with
the definition \eqref{eq:definvar} of invariance of a family $I$,
yields
\begin{linenomath} 
\begin{equation*}
 f_{I}(\bi{y}|\bs{\lambda})=
 f_I(\bi{g}_{a}^{-1}(\bi{y})|
 \bar{\bi{g}}_{a}^{-1}(\bs{\lambda}))\,
 \mbox{\large $|$}\partial_{\bi{y}}\bi{g}_{a}^{-1}(\bi{y})\mbox{\large $|$}
\end{equation*} 
\end{linenomath}
for all $\bi{y}\equiv \bi{g}_{a}(\bi{x})$ such that
$\bi{x}\in V_{\bi{X}}(\bs{\theta})$, where 
$\bs{\lambda}\equiv\bar{\bi{g}}_{a}(\bs{\theta})$ and
  $\bar{\bi{g}}_{a}\in\bar{\mathcal G}$. 
%

\section{Inverse probability distributions}
\label{sec:inverse}

\begin{define}[Inverse probability distributions]
\label{def:inversepr}
Suppose there exist probability spaces 
$\left(\Omega_{\bs{\theta}},\Sigma_{\bs{\theta}},P\right)$, 
$\Omega_{\bs{\theta}}\subset\Omega$, for all $\bs{\theta}\in
V_{\mathbf{\Theta}}$  and a random variable 
$(\mathbf{\Theta},\bi{X}):\Omega_{\bs{\theta}}\longrightarrow\left(\bs{\theta},
\mathbb{R}^n\right)$ that together
lead to the parametric family $I$ of continuous 
direct probability distributions
$Pr_I\left(S_{\mathbf{\Theta}=\bs{\theta},\bi{X}}|\bs{\theta}\right)$, 
$S_{\mathbf{\Theta}=\bs{\theta},\bi{X}}\in{\mathcal B}^n$ whose pdf's
are denoted by $f_I(\bi{x}|\bs{\theta})$.
Let, in addition, for some of those realizations $\bi{x}$ of $\bi{X}$ 
for which 
\begin{linenomath} 
\begin{equation}
 \int_{V_{\mathbf{\Theta}}}f_I(\bi{x}|\bs{\theta})\,d^m\bs{\theta}
 > 0 \ ,
\label{eq:condexist}
\end{equation} 
\end{linenomath}
there exist also probability spaces 
$\left(\Omega_{\bi{x}},\Sigma_{\bi{x}},P\right)$, 
$\Omega_{\bi{x}}\subset\Omega$, such that the function 
$(\mathbf{\Theta},\bi{X}):
\Omega_{\bi{x}}\longrightarrow\left(V_{\mathbf{\Theta}},\bi{x}\right)$
is $\Sigma_{\bi{x}}$-measurable (i.e., 
$A_{\mathbf{\Theta}\le\bs{\theta}}\equiv 
\{\omega\in\Omega_{\bi{x}} :
\mathbf{\Theta}\le\bs{\theta}\} \in \Sigma_{\bi{x}}$ for all
$\bs{\theta}\in V_{\mathbf{\Theta}}$) and thus 
a random variable also on
$\left(\Omega_{\bi{x}},\Sigma_{\bi{x}},P\right)$. Then, the
probability distributions, resulting from the probability
spaces $\left(\Omega_{\bi{x}},\Sigma_{\bi{x}},P\right)$ and from 
the corresponding random variable $(\mathbf{\Theta},\bi{X})$, 
are called {\em inverse
  probability distributions}. The cdf's and the pdf's
that correspond to the inverse probability
distributions
are denoted by $F_I\left(\bs{\theta}|\bi{x}\right)$ and 
$f_I\left(\bs{\theta}|\bi{x}\right)$, respectively.

Likewise, let $(\mathbf{\Theta},\bi{X})$ be further partitioned into
 $(\mathbf{\Theta}_1,\mathbf{\Theta}_2,\bi{X})$ and let for some of
those realizations $\bs{\theta}_1$ and $\bi{x}$ for which
\begin{linenomath} 
\begin{equation}
 \int_{V_{\mathbf{\Theta}_1,\bs{\theta}_2}}
  f_I(\bi{x}|\bs{\theta}_1,\bs{\theta}_2)\,d^{m_1}\bs{\theta}_1
 > 0 \ ,
\label{eq:condexist1}
\end{equation} 
\end{linenomath}
there exist probability spaces 
$\left(\Omega_{\bs{\theta}_2,\bi{x}},\Sigma_{\bs{\theta}_2,\bi{x}},P\right)$
such that the function $(\mathbf{\Theta}_1,\mathbf{\Theta}_2,\bi{X}):
\Omega_{\bs{\theta}_2,\bi{x}}\longrightarrow
\left(V_{\mathbf{\Theta}_1,\bs{\theta}_2},\bs{\theta}_2,\bi{x}\right)$
is $\Sigma_{\bs{\theta}_2,\bi{x}}$-measurable, 
$A_{\mathbf{\Theta}_1\le\bs{\theta}_1}\equiv 
\{\omega\in\Omega_{\bs{\theta}_1,\bi{x}} :
\mathbf{\Theta}_1\le$ $\bs{\theta}_1\} \in
\Sigma_{\bs{\theta}_2,\bi{x}}$  
for all $(\bs{\theta}_1,\bs{\theta}_2)\in V_{\mathbf{\Theta}_1,\bs{\theta}_2}$.
Then, the cdf's and the pdf's that correspond to the resulting inverse
probability distributions are denoted by 
$F_I\left(\bs{\theta}_1|\bs{\theta}_2,\bi{x}\right)$ and 
$f_I\left(\bs{\theta}_1|\bs{\theta}_2,\bi{x}\right)$, respectively.
\end{define}
\begin{remark}
\label{rem:non-vanint}
{\rm The integrals \eqref{eq:condexist} and \eqref{eq:condexist1} need 
not be finite. The reasons for requiring the two integrals 
to be strictly positive will become 
apparent within the context of Proposition\;\ref{theo:consistency}, below.
}
\end{remark}

Apart from the direct and the inverse probability distributions,
their mixtures may also exist. For example,
$F_I\left(\bs{\theta},\bi{x}_1|\bi{x}_2\right)$,  
$F_I\left(\bs{\theta}_1,\bi{x}_1|\bs{\theta}_2,\bi{x}_2\right)$, 
$f_I\left(\bs{\theta},\bi{x}_1|\bi{x}_2\right)$ 
and $f_I\left(\bs{\theta}_1,\bi{x}_1|\bs{\theta}_2,\bi{x}_2\right)$ 
are the cdf's and
the pdf's of two of the distributions that are neither purely direct nor
purely inverse.

From a mathematical perspective, the direct and the inverse
probability distributions, as well as their mixtures, 
share identical properties, some of which
were discussed in Section\;\ref{ss:general}.
The following three rules that apply to inverse probability
distributions are
obtained by invoking the equivalence between the two types 
of distributions.
\begin{rul}[Parameter transformation]
Let $f_I(\bs{\theta}|\bi{x})$ be a pdf of an inverse probability
distribution and let $(\bar{\bi{s}},\bi{s}):
(\mathbf{\Theta},\bi{X})\longrightarrow
(\bar{\bi{s}}\circ\mathbf{\Theta},\bi{s}\circ\bi{X})\equiv
(\mathbf{\Lambda},\bi{Y})$ be a differentiable transformation with a
non-vanishing Jacobian on the entire support of
$f_I(\bs{\theta}|\bi{x})$. Then, an inverse pdf 
$f_{I'}(\bs{\lambda}|\bi{y})$ also exists and is related to
$f_I(\bs{\theta}|\bi{x})$ as
\begin{linenomath} 
\begin{equation}
 f_{I'}(\bs{\lambda}|\bi{y})
=
f_{I}(\bar{\bi{s}}^{-1}(\bs{\lambda})|\bi{s}^{-1}(\bi{y}))\,
 \mbox{\large $|$}\partial_{\bs{\lambda}}\bar{\bi{s}}^{-1}(\bs{\lambda})\mbox{\large $|$}  \ .
\label{eq:condtrans2}
\end{equation} 
\end{linenomath}
Similarly, when there exist an inverse pdf 
$f_I(\bs{\theta}_1|\bs{\theta}_2,\bi{x})$ and a differentiable
transformation $(\bar{\bi{s}}_1,\bar{\bi{s}}_2,\bi{s}):
(\mathbf{\Theta}_1,\mathbf{\Theta}_2,\bi{X})\longrightarrow
(\bar{\bi{s}}_1\circ\mathbf{\Theta},\bar{\bi{s}}_2\circ\mathbf{\Theta},
\bi{s}\circ\bi{X})\equiv
(\mathbf{\Lambda}_1,\mathbf{\Lambda}_2,\bi{Y})$ with a non-vanishing
Jacobian on the support of $f_I(\bs{\theta}_1|\bs{\theta}_2,\bi{x})$,
there exists a pdf $f_{I'}(\bs{\lambda}_1|\bs{\lambda}_2,\bi{x})$ such that
\begin{linenomath} 
\begin{equation}
\begin{split}
&f_{I'}(\bs{\lambda}_1|\bs{\lambda}_2,\bi{y})
= \\
&f_{I}(\bar{\bi{s}}_1^{-1}(\bs{\lambda}_1)|
 \bar{\bi{s}}_2^{-1}(\bs{\lambda}_2),\bi{s}^{-1}(\bi{y}))\,
 \mbox{\large $|$}\partial_{\bs{\lambda}_1}\bar{\bi{s}}_1^{-1}(\bs{\lambda}_1)\mbox{\large $|$}  \ .
\end{split}
\label{eq:condtrans3}
\end{equation} 
\end{linenomath}
\end{rul}
\hskip 0mm \textit{Proof.} \;If $f_I(\bs{\theta},\bi{x})$ exists, 
   equation \eqref{eq:condtrans2} follows from \eqref{eq:condtrans0}
   by substitutions $\bi{X}_1\to\mathbf{\Theta}$, $\bi{X}_2\to\bi{X}$,
   $\bi{Y}_1\to\mathbf{\Lambda}$, $\bi{Y}_2\to\bi{Y}$,
   $\bi{s}_1\to\bar{\bi{s}}$ and $\bi{s}_2\to\bi{s}$.
   Similarly, if $f_I(\bs{\theta}_1,\bs{\theta}_2,\bi{x})$ exists, 
   \eqref{eq:condtrans3} is deduced from
   \eqref{eq:condtrans0}
   by substitutions $\bi{X}_1\to\mathbf{\Theta}_1$,
   $\bi{X}_2\to(\mathbf{\Theta}_2,\bi{X})$, 
   $\bi{Y}_1\to\mathbf{\Lambda}_1$, $\bi{Y}_2\to(\mathbf{\Lambda}_2,\bi{Y})$,
   $\bi{s}_1\to\bar{\bi{s}}_1$ and
   $\bi{s}_2\to(\bar{\bi{s}}_2,\bi{s})$. If, on the other hand, 
   the joint pdf's $f_I(\bs{\theta},\bi{x})$ and
   $f_I(\bs{\theta}_1,\bs{\theta}_2,\bi{x})$ do not exist, 
   equations  \eqref{eq:condtrans2} and  \eqref{eq:condtrans3} are
   definitions for $f_{I'}(\bs{\lambda}|\bi{y})$ and 
   $f_{I'}(\bs{\lambda}_1|\bs{\lambda}_2,\bi{y})$, respectively, in the
   same way as $f_{\mathbf{Y}}(\bi{y}_1|\bi{y}_2)$ was defined by
   \eqref{eq:condtrans1}. 
%
\hfill $\Box$
\vskip 3mm 
\begin{rul}[Product rule]
\label{prop:prod}
Let there exist an inverse pdf
$f_I(\bs{\theta}_1,\bs{\theta}_2|\bi{x})$ 
and the corresponding mar\-gi\-nal pdf 
\begin{linenomath} 
\begin{equation*}
 f_I(\bs{\theta}_2|\bi{x})\equiv\int_{V_{\mathbf{\Theta}_1}}
 f_I(\bs{\theta}_1,\bs{\theta}_2|\bi{x})\,d^{m_1}\bs{\theta}_1 \ .
\end{equation*} 
\end{linenomath}
Then, for all 
$\bs{\theta}_2$ and $\bi{x}$ for which $f_I(\bs{\theta}_2|\bi{x})>0$,
\begin{linenomath} 
\begin{equation}
 f_I(\bs{\theta}_1|\bs{\theta}_2,\bi{x})=
 \frac{f_I(\bs{\theta}_1,\bs{\theta}_2|\bi{x})}
      {f_I(\bs{\theta}_2|\bi{x})} 
\label{eq:invproduct}
\end{equation} 
\end{linenomath}
holds uniquely (Lebesgue measure) 
$\nu_L$-almost everywhere on $\mathbb{R}^{m_1}$.
\end{rul}
\hskip 0mm \textit{Proof.} \;The product rule \eqref{eq:invproduct} follows
   immediately from \eqref{eq:condrat1.a} by making substitutions
   $\bi{x}_1\to\bs{\theta}_1$, $\bi{x}_2\to\bs{\theta}_2$ and
   $\bi{x}_3\to\bi{x}$.
%
\hfill $\Box$ 
%
%
\begin{rul}[Bayes' Theorem]
Let a random vector be partitioned into
  $(\mathbf{\Theta}_1,\mathbf{\Theta}_2,\bi{X}_1,$ $\bi{X}_2)$, let
there exist pdf's
$f_I\left(\bs{\theta}_1,\bi{x}_1|\bs{\theta_2},\bi{x}_2\right)$ and
$f_I\left(\bs{\theta}_1,\bi{x}_2|\bs{\theta}_2,\bi{x}_1\right)$, 
let marginal pdf's $f_I\left(\bi{x}_1|\bs{\theta}_2,\bi{x}_2\right)$,
$f_I\left(\bi{x}_2|\bs{\theta}_2,\bi{x}_1\right)$,
$f_I\left(\bs{\theta}_1|\bs{\theta}_2,\bi{x}_2\right)$ and
$f_I\left(\bs{\theta}_1|\bs{\theta}_2,\bi{x}_1\right)$ be
non-vanishing, 
and let the components $\bi{X}_1$ and $\bi{X}_2$ of the partition
be independent random variables:
$f_I\left(\bi{x}_1|\bs{\theta}_1,\bs{\theta}_2,\bi{x}_2\right)=
f_I\left(\bi{x}_1|\bs{\theta}_1,\bs{\theta}_2\right)$ and
$f_I\left(\bi{x}_2|\bs{\theta}_1,\bs{\theta}_2,\bi{x}_1\right)=
f_I\left(\bi{x}_2|\bs{\theta}_1,\bs{\theta}_2\right)$.
Then, there exists a conditional pdf
$f_I\left(\bs{\theta}_1|\bs{\theta}_2,\bi{x}_1,\bi{x}_2\right)$ such that 
\begin{linenomath} 
\begin{equation}
\begin{split}
 f_I\left(\bs{\theta}_1|\bs{\theta}_2,\bi{x}_1,\bi{x}_2\right)&=
 \frac{f_I\left(\bs{\theta}_1|\bs{\theta}_2,\bi{x}_2\right)\,
 f_I\left(\bi{x}_1|\bs{\theta}_1,\bs{\theta}_2\right)}
      {f_I\left(\bi{x}_1|\bs{\theta}_2,\bi{x}_2\right)} \\
&=
 \frac{f_I\left(\bs{\theta}_1|\bs{\theta}_2,\bi{x}_1\right)\,
       f_I\left(\bi{x}_2|\bs{\theta}_1,\bs{\theta}_2\right)}
      {f_I\left(\bi{x}_2|\bs{\theta}_2,\bi{x}_1\right)} \ .
\end{split}
\label{eq:bayes1}
\end{equation} 
\end{linenomath}
If, on the other hand, a random vector is partitioned into 
$(\mathbf{\Theta},\bi{X}_1,\bi{X}_2)$,
\begin{linenomath} 
\begin{equation}
\begin{split}
 f_I\left(\bs{\theta}|\bi{x}_1,\bi{x}_2\right)&=
 \frac{f_I\left(\bs{\theta}|\bi{x}_2\right)\,
       f_I\left(\bi{x}_1|\bs{\theta}\right)}
      {f_I\left(\bi{x}_1|\bi{x}_2\right)} \\
&=
 \frac{f_I\left(\bs{\theta}|\bi{x}_1\right)\,
       f_I\left(\bi{x}_2|\bs{\theta}\right)}
      {f_I\left(\bi{x}_2|\bi{x}_1\right)} 
\end{split}
\label{eq:bayes2}
\end{equation} 
\end{linenomath}
holds true under analogous conditions.
\end{rul}
\hskip 0mm \textit{Proof.} \;Equation \eqref{eq:bayes1} follows
   from \eqref{eq:bayes0} by making substitutions
   $\bi{y}\to\bs{\theta}_1$, $\bi{z}\to\bs{\theta}_2$,
   $\bi{t}\to\bi{x}_1$ and $\bi{w}\to\bi{x}_2$,
   whereas \eqref{eq:bayes2} is obtained from \eqref{eq:bayes01} by
   substitutions $\bi{y}\to\bs{\theta}$, $\bi{t}\to\bi{x}_1$ 
   and $\bi{w}\to\bi{x}_2$.
\hfill $\Box$ 
\vskip 3mm \noindent
Equations \eqref{eq:bayes1} and \eqref{eq:bayes2} are also
referred to as Bayes' Theorem (\cite{bay,lap}) or the principle of 
inverse probability (\cite{jef}, \S\,1.22, p.\,28),
written in terms of pdf's. In the equations, 
$f_I\left(\bs{\theta}_1|\bs{\theta}_2,\bi{x}_1,\bi{x}_2\right)$
and $f_I\left(\bs{\theta}|\bi{x}_1,\bi{x}_2\right)$ are
called the {\em posterior pdf's}, 
$f_I\left(\bs{\theta}_1|\bs{\theta}_2,\bi{x}_{1,2}\right)$
and
$f_I\left(\bs{\theta}|\bi{x}_{1,2}\right)$ are the 
so-called {\em prior pdf's}, 
$f_I\left(\bi{x}_{1,2}|\bs{\theta}_1,\bs{\theta}_2\right)$ 
and $f_I\left(\bi{x}_{1,2}|\bs{\theta}\right)$ are  
the {\em likelihood densities}, while 
$f_I\left(\bi{x}_{1,2}|\bs{\theta}_2,\bi{x}_{2,1}\right)$
and $f_I\left(\bi{x}_{1,2}|\bi{x}_{2,1}\right)$ are 
the {\em predictive pdf's}.  While the predictive pdf's 
are determined by the normalization condition on the posterior pdf's, 
e.g.,
\begin{linenomath} 
\begin{equation*}
\begin{split}
&f_I\left(\bi{x}_{1,2}|\bs{\theta}_2,\bi{x}_{2,1}\right)= \\
&\int_{V_{\mathbf{\Theta}_1}}
 f_I\left(\bs{\theta}_1|\bs{\theta}_2,\bi{x}_{2,1}\right)\,
 f_I\left(\bi{x}_{1,2}|\bs{\theta}_1,\bs{\theta}_2\right)\,d^{m_1}\bs{\theta}_1 
\ ,
\end{split}
\end{equation*} 
\end{linenomath}
the general form of the prior pdf's 
$f_I\left(\bs{\theta}_1|\bs{\theta}_2,\bi{x}_{1,2}\right)$
and $f_I\left(\bs{\theta}|\bi{x}_{1,2}\right)$
is prescribed by the following Proposition.
\begin{theo}
\label{theo:consistency}
Suppose that conditions for Bayes' Theorem \eqref{eq:bayes1} 
are fulfilled: a random vector is partitioned into
  $(\mathbf{\Theta}_1,\mathbf{\Theta}_2,\bi{X}_1,$ $\bi{X}_2)$,
there exist conditional pdf's 
$f_I\left(\bs{\theta}_1,\bi{x}_1|\bs{\theta_2},\bi{x}_2\right)$ and
$f_I\left(\bs{\theta}_1,\bi{x}_2|\bs{\theta}_2,\bi{x}_1\right)$, 
the marginal pdf's $f_I\left(\bi{x}_1|\bs{\theta}_2,\bi{x}_2\right)$,
$f_I\left(\bi{x}_2|\bs{\theta}_2,\bi{x}_1\right)$,
$f_I\left(\bs{\theta}_1|\bs{\theta}_2,\bi{x}_2\right)$ and
$f_I\left(\bs{\theta}_1|\bs{\theta}_2,\bi{x}_1\right)$ are
positive, 
and the components $\bi{X}_1$ and $\bi{X}_2$ of the partition
are independent random variables with identical probability
distribution. In addition, let $V_{\mathbf{\Theta}}=
(V_{\mathbf{\Theta}_1},V_{\mathbf{\Theta}_2})$ stand for the space 
of the parameter $\mathbf{\Theta}=(\mathbf{\Theta}_1,
\mathbf{\Theta}_2)$ and let
$\widetilde{V}_{\mathbf{\Theta}_1}(\bi{x}_{1,2},\bs{\theta}_2)\equiv
\left\{\bs{\theta}_1\in V_{\mathbf{\Theta}_1}:
f_I(\bi{x}_{1,2}|\bs{\theta}_1,\bs{\theta}_2) > 0\right\}$.
Then, for $\bs{\theta}_1\in 
\widetilde{V}_{\mathbf{\Theta}_1}(\bi{x}_{1,2},\bs{\theta}_2)$,
\begin{linenomath} 
\begin{equation}
f_I(\bs{\theta}_1|\bs{\theta}_2,\bi{x}_{1,2})=
\frac{\zeta_{I,\mathbf{\Theta}_1|\bs{\theta}_2}(\bs{\theta}_1,\bs{\theta}_2)}
{\eta_{I,\mathbf{\Theta}_1|\bs{\theta}_2}(\bi{x}_{1,2},\bs{\theta}_2,)}
\,f_I\left(\bi{x}_{1,2}|\bs{\theta}_1,\bs{\theta}_2\right) 
\label{eq:bayprime1}
\end{equation} 
\end{linenomath}
is the most general form of the pdf's 
$f_I\left(\bs{\theta}_1|\bs{\theta}_2,\bi{x}_{1,2}\right)$. Similarly,
when a random vector is partitioned into 
$(\mathbf{\Theta},\bi{X}_1,$ $\bi{X}_2)$, the conditions 
for Bayes' Theorem \eqref{eq:bayes2} are fulfilled and  
$\bs{\theta}\in\widetilde{V}_{\mathbf{\Theta}}(\bi{x}_{1,2})$,
 $\widetilde{V}_{\mathbf{\Theta}}(\bi{x}_{1,2})\equiv
\left\{\bs{\theta}\in V_{\mathbf{\Theta}}:
f_I(\bi{x}_{1,2}|\bs{\theta}) > 0\right\}$,
\begin{linenomath} 
\begin{equation}
f_I\left(\bs{\theta}|\bi{x}_{1,2}\right)=
\frac{\zeta_{I,\mathbf{\Theta}}(\bs{\theta})}
{\eta_{I,\mathbf{\Theta}}(\bi{x}_{1,2})}
\,f_I\left(\bi{x}_{1,2}|\bs{\theta}\right) 
\label{eq:bayprime2}
\end{equation} 
\end{linenomath}
is the most general form of the pdf's 
$f_I\left(\bs{\theta}|\bi{x}_{1,2}\right)$.
The functions
$\zeta_{I,\mathbf{\Theta}_1|\bs{\theta}_2}(\bs{\theta}_1,\bs{\theta}_2)$
and $\zeta_{I,\mathbf{\Theta}}(\bs{\theta})$
in equations \eqref{eq:bayprime1} and \eqref{eq:bayprime2} 
are called the {\em consistency factors}.
\end{theo}
Domains of
$f_I\left(\bs{\theta}_1|\bs{\theta}_2,\bi{x}_{1,2}\right)$ 
and $f_I\left(\bs{\theta}|\bi{x}_{1,2}\right)$ are extended beyond
the supports $V_{\bi{X}_1}(\bs{\theta}_1,\bs{\theta}_2) = 
V_{\bi{X}_2}(\bs{\theta}_1,\bs{\theta}_2)$ on which
$f_I\left(\bi{x}_{1,2}|\bs{\theta}_1,\bs{\theta}_2\right)$ and
$f_I\left(\bi{x}_{1,2}|\bs{\theta}\right)$ are positive by defining
\begin{linenomath} 
\begin{equation*}
f_I\left(\bs{\theta}_1|\bs{\theta}_2,\bi{x}_{1,2}\right)\equiv
\frac{\zeta_{I,\mathbf{\Theta}_1|\bs{\theta}_2}(\bs{\theta}_1,\bs{\theta}_2)}
{\eta_{I,\mathbf{\Theta}_1|\bs{\theta}_2}(\bi{x}_{1,2},\bs{\theta}_2,)}
\,f_I\left(\bi{x}_{1,2}|\bs{\theta}_1,\bs{\theta}_2\right) 
\end{equation*} 
\end{linenomath}
for all $\bi{x}_{1,2}\not\in
V_{\bi{X}_{1,2}}(\bs{\theta}_1,\bs{\theta}_2)$ and 
\begin{linenomath} 
\begin{equation*}
f_I(\bs{\theta}|\bi{x}_{1,2})\equiv
\frac{\zeta_{I,\mathbf{\Theta}}(\bs{\theta})}
{\eta_{I,\mathbf{\Theta}}(\bi{x}_{1,2})}
\,f_I\left(\bi{x}_{1,2}|\bs{\theta}\right) 
\end{equation*} 
\end{linenomath}
for all $\bi{x}_{1,2}\not\in V_{\bi{X}_{1,2}}(\bs{\theta})$. For the
sake of symmetry between the direct and the inverse probability
distributions, the domains of the inverse pdf's may be extended
even further by defining 
$f_I\left(\bs{\theta}_1|\bs{\theta}_2,\bi{x}_{1,2}\right)\equiv 0$
for $(\bs{\theta}_1,\bs{\theta}_2)\notin V_{\mathbf{\Theta}_1,\bs{\theta}_2}$
and $f_I\left(\bs{\theta}|\bi{x}_{1,2}\right)\equiv 0$
for $\bs{\theta}\notin V_{\mathbf{\Theta}}$. In this way, the inverse
probability distribution spaces $(V_{\bi{x}},\Sigma_{\bi{x}},Pr)$ and
$(V_{\bs{\theta}_2,\bi{x}},\Sigma_{\bs{\theta}_2,\bi{x}},Pr)$ are also 
extended to $(\mathbb{R}^m,{\mathcal B}^m,Pr)$ and
 $(\mathbb{R}^{m_1},{\mathcal B}^{m_1},Pr)$, respectively.
Then, the {\em normalization factors} 
$\eta_{I,\mathbf{\Theta}_1|\bs{\theta}_2}(\bi{x}_{1,2},\bs{\theta}_2,)$ 
and $\eta_{I,\mathbf{\Theta}}(\bi{x}_{1,2})$ are
determined by invoking normalization of the 
pdf's $f_I\left(\bs{\theta}_1|\bs{\theta}_2,\bi{x}_{1,2}\right)$ 
and $f_I\left(\bs{\theta}|\bi{x}_{1,2}\right)$:
\begin{linenomath} 
\begin{equation*}
\begin{split}
&\eta_{I,\mathbf{\Theta}_1|\bs{\theta}_2}(\bi{x}_{1,2},\bs{\theta}_2,)
=\\
&\int_{\mathbb{R}^{m_1}}
 \zeta_{I,\mathbf{\Theta}_1|\bs{\theta}_2}(\bs{\theta}_1,\bs{\theta}_2)
 \,f_I\left(\bi{x}_{1,2}|\bs{\theta}_1,\bs{\theta}_2\right)\,
 d^{m_1}\bs{\theta}_1
\end{split}
\end{equation*} 
\end{linenomath}
and
\begin{linenomath} 
\begin{equation*}
\eta_{I,\mathbf{\Theta}}(\bi{x}_{1,2})
=
\int_{\mathbb{R}^m}\zeta_{I,\mathbf{\Theta}}(\bs{\theta})
 \,f_I\left(\bi{x}_{1,2}|\bs{\theta}\right)
 d^{m}\bs{\theta} \ .
\end{equation*} 
\end{linenomath}
Non-vanishing integrals \eqref{eq:condexist} and \eqref{eq:condexist1}
thus represent necessary
conditions for normalizability \eqref{eq:norpdf} of the inverse pdf's
$f_I\left(\bs{\theta}|\bi{x}_{1,2}\right)$ and 
$f_I\left(\bs{\theta}_1|\bs{\theta}_2,\bi{x}_{1,2}\right)$.

For discrete random variables $\bi{X}_1$ and $\bi{X}_2$,
the appropriate forms of the pdf's \linebreak
$f_I\left(\bs{\theta}_1|\bs{\theta}_2,\bi{x}_{1,2}\right)$ and 
$f_I\left(\bs{\theta}|\bi{x}_{1,2}\right)$ 
are obtained 
by replacing the likelihood densities
$f_I\left(\bi{x}_{1,2}|\bs{\theta}_1,\bs{\theta}_2\right)$
and $f_I\left(\bi{x}_{1,2}|\bs{\theta}\right)$ in \eqref{eq:bayprime1}
and \eqref{eq:bayprime2}
with the {\em probability mass functions} 
$p_I\left(\bi{x}_{1,2}|\bs{\theta}_1,\bs{\theta}_2\right)$
and $p_I\left(\bi{x}_{1,2}|\bs{\theta}\right)$ that coincide
with probability distributions for the points
$\bi{X}_{1,2}=\bi{x}_{1,2}$ of a state space $\mathbb{R}^{n}$
of the variables $\bi{X}_1$ and $\bi{X}_2$, given the realizations
$(\mathbf{\Theta}_1,\mathbf{\Theta}_2)=(\bs{\theta}_1,\bs{\theta}_2)$ 
and $\mathbf{\Theta}=\bs{\theta}$ of the corresponding parameters.
\begin{remark}
\label{rem:bayesvscons}
{\rm In equations \eqref{eq:bayprime1} and \eqref{eq:bayprime2}, the pdf's 
     $f_I\left(\bs{\theta}_1|\bs{\theta}_2,\bi{x}_{1,2}\right)$ and
     $f_I\left(\bs{\theta}|\bi{x}_{1,2}\right)$ 
     are directly proportional to the pdf's
     $f_I\left(\bi{x}_{1,2}|\bs{\theta}_1,\bs{\theta}_2\right)$
     and $f_I\left(\bi{x}_{1,2}|\bs{\theta}\right)$ of the
     corresponding direct probability distributions. This is very 
     similar to equations \eqref{eq:bayes1} and \eqref{eq:bayes2}
     of Bayes' Theorem with the posterior pdf's
     $f_I\left(\bs{\theta}_1|\bs{\theta}_2,\bi{x}_1,\bi{x}_2\right)$
     and $f_I\left(\bs{\theta}|\bi{x}_1,\bi{x}_2\right)$ being
     proportional to the likelihood densities 
     $f_I\left(\bi{x}_{1,2}|\bs{\theta}_1,\bs{\theta}_2\right)$ 
     and $f_I\left(\bi{x}_{1,2}|\bs{\theta}\right)$. 
     But there is also a fundamental difference between the equations
     of Bayes' Theorem and those of Proposition\;\ref{theo:consistency}: 
     while the proportionality coefficients 
     $f_I\left(\bs{\theta}_1|\bs{\theta}_2, \bi{x}_{1,2}\right)$ and
     $f_I\left(\bs{\theta}|\bi{x}_{1,2}\right)$ between the
     posterior pdf's and the likelihood densities in Bayes' Theorem
     are the prior pdf's, the consistency factors 
     $\zeta_{I,\mathbf{\Theta}_1|\bs{\theta}_2}(\bs{\theta}_1,\bs{\theta}_2)$
     and $\zeta_{I,\mathbf{\Theta}}(\bs{\theta})$ that are 
     proportionality coefficients between the inverse and the 
     direct pdf's in \eqref{eq:bayprime1} and \eqref{eq:bayprime2}
      need not be congruent with 
     all the properties of probability density functions and should therefore
     not be confused with the so-called {\em non-informative prior pdf's}
     $f_I(\bs{\theta}_{1,2}|\bs{\theta}_{2,1})$ and 
     $f_I(\bs{\theta})$ (see also Section\;\ref{ss:unique}, below). 
     The properties of the consistency factors 
     are extensively discussed in the next section. 
}
\end{remark}
%

\section{The consistency factors}
\label{sec:factors}


\subsection{General properties of the consistency factors}
\label{ss:properties}
According to Proposition\;\ref{theo:consistency}, for a consistent
assignment of inverse probability distributions, the appropriate
consistency factors $\zeta_{I,\mathbf{\Theta}}(\bs{\theta})$ and 
$\zeta_{I,\mathbf{\Theta}_{1,2}|\bs{\theta}_{2,1}}(\bs{\theta}_1,\bs{\theta}_2)$
need be uniquely determined. In what follows, we discuss some
of the properties of the consistency factors that will be invoked
during their determination.
\begin{proper}[Uniqueness]
\label{prop:unique}
A consistency factors $\zeta_{I,\mathbf{\Theta}}(\bs{\theta})$ can only be
determined up to a factor $\chi_{I,\mathbf{\Theta}}(\bi{x}_{1,2})$
that is an arbitrary function of $\bi{x}_{1,2}$. Also,
$\zeta_{I,\mathbf{\Theta}_{1,2}|\bs{\theta}_{2,1}}(\bs{\theta}_1,\bs{\theta}_2)$
is determined only up to an arbitrary multiplier 
$\chi_{I,\mathbf{\Theta}_{1,2}|\bs{\theta}_{2,1}}(\bi{x}_{1,2},\bs{\theta}_{2,1})$.
\end{proper}
\hskip 0mm \textit{Proof.} \;Multiplying
$\zeta_{I,\mathbf{\Theta}}(\bs{\theta})$ by 
$\chi_{I,\mathbf{\Theta}}(\bi{x}_{1,2})$ results in
multiplying $\eta_{I,\mathbf{\Theta}}(\bi{x}_{1,2})$ by the same
factor,
such that the factor cancels in the ratio 
$\zeta_{I,\mathbf{\Theta}}(\bs{\theta})/\eta_{I,\mathbf{\Theta}}(\bi{x}_{1,2})$.
Identical arguments apply when 
$\zeta_{I,\mathbf{\Theta}_{1,2}|\bs{\theta}_{2,1}}(\bs{\theta}_1,\bs{\theta}_2)$
is multiplied by 
$\chi_{I,\mathbf{\Theta}_{1,2}|\bs{\theta}_{2,1}}(\bi{x}_{1,2},\bs{\theta}_{2,1})$.
%
\hfill $\Box$
\vskip 3mm 
\begin{proper}[Sign]
\label{prop:sign}
A consistency factor $\zeta_{I,\mathbf{\Theta}}(\bs{\theta})$ is either
positive or negative on the parameter space $V_{\mathbf{\Theta}}$,
and so is
$\zeta_{I,\mathbf{\Theta}_{1,2}|\bs{\theta}_{2,1}}(\bs{\theta}_1,\bs{\theta}_2)$
on $V_{\mathbf{\Theta}_{1,2}|\bs{\theta}_{2,1}}$.
\end{proper}
\hskip 0mm \textit{Proof.} \;The normalization factors 
$\eta_{I,\mathbf{\Theta}}(\bi{x}_{1,2})$ are either positive or negative,
and the pdf's $f_I\left(\bs{\theta}|\bi{x}_{1,2}\right)$ and 
$f_I\left(\bi{x}_{1,2}|\bs{\theta}\right)$ are non-negative, such
that $\zeta_{I,\mathbf{\Theta}}(\bs{\theta})$ must be of the same sign as 
$\eta_{I,\mathbf{\Theta}}(\bi{x}_{1,2})$, i.e., either positive or
negative for all $\bs{\theta}\in V_{\mathbf{\Theta}}$. The 
same holds true for 
$\eta_{I,\mathbf{\Theta}_{1,2}|\bs{\theta}_{2,1}}(\bi{x}_{1,2},\bs{\theta}_{2,1})$,
 $f_I\left(\bs{\theta}_{1,2}|\bs{\theta}_{2,1},\bi{x}_{1,2}\right)$,
$f_I\left(\bi{x}_{1,2}|\bs{\theta}_1,\bs{\theta}_2\right)$ and
$\zeta_{I,\mathbf{\Theta}_{1,2}|\bs{\theta}_{2,1}}(\bs{\theta}_1,\bs{\theta}_2)$.
\hfill $\Box$
\vskip 3mm 
\begin{proper}[Transformations]
\label{prop:trans}
Suppose that the premises of Proposition\;\ref{theo:consistency} are fulfilled
such that pdf's $f_I\left(\bs{\theta}|\bi{x}_{1,2}\right)$ and 
$f_I\left(\bi{x}_{1,2}|\bs{\theta}\right)$ are related according to 
\eqref{eq:bayprime2}. Let, in addition,  $(\bar{\bi{s}},\bi{s}):
(\mathbf{\Theta},\bi{X})\longrightarrow
(\bar{\bi{s}}\circ\mathbf{\Theta},\bi{s}\circ\bi{X})\equiv
(\mathbf{\Lambda},\bi{Y})$ be a differentiable transformation with
non-vanishing Jacobians
$\mbox{\large $|$}\partial_{\bs{\lambda}}\bar{\bi{s}}(\bs{\lambda})\mbox{\large $|$}$  
and $\mbox{\large $|$}
\partial_{\bi{x}_{1,2}}\bi{s}(\bi{x}_{1,2})\mbox{\large $|$}$ for all $\bs{\theta}$
and $\bi{x}_{1,2}$ for which 
$f_I\left(\bi{x}_{1,2}|\bs{\theta}\right)$
is positive. Then, the consistency and the normalization factors
that relate $f_{I'}\left(\bs{\lambda}|\bi{y}_{1,2}\right)$ and 
$f_{I'}\left(\bi{y}_{1,2}|\bs{\lambda}\right)$ read
\begin{linenomath} 
\begin{equation}
\zeta_{I',\mathbf{\Lambda}}(\bs{\lambda})= \chi_{I',\mathbf{\Lambda}}\,
\zeta_{I,\mathbf{\Theta}}[\bar{\bi{s}}^{-1}(\bs{\lambda})]\,
 \mbox{\large $|$}\partial_{\bs{\lambda}}\bar{\bi{s}}^{-1}(\bs{\lambda})\mbox{\large $|$}
\label{eq:cftr1}
\end{equation} 
\end{linenomath}
and
\begin{linenomath} 
\begin{equation}
\eta_{I',\mathbf{\Lambda}}(\bi{y}_{1,2})= \chi_{I',\mathbf{\Lambda}}\,
\eta_{I,\mathbf{\Theta}}[\bi{s}^{-1}(\bi{y}_{1,2})]\,
 \mbox{\large $|$}\partial_{\bi{y}_{1,2}}\bi{s}^{-1}(\bi{y}_{1,2})\mbox{\large $|$} \ .
\label{eq:nftr1}
\end{equation} 
\end{linenomath}
Similarly, 
for $f_I\left(\bs{\theta}_{1,2}|\bs{\theta}_{2,1},\bi{x}_{1,2}\right)$ and 
$f_I\left(\bi{x}_{1,2}|\bs{\theta}_1,\bs{\theta}_2\right)$, 
\begin{linenomath} 
\begin{equation}
\begin{split}
&\zeta_{I',\mathbf{\Lambda}_{1,2}|\bs{\lambda}_{2,1}}(\bs{\lambda}_1,\bs{\lambda}_{2})
= \chi_{I',\mathbf{\Lambda}_{1,2}|\bs{\lambda}_{2,1}}\,\times \\
&\zeta_{I,\mathbf{\Theta}_{1,2}|\bar{\bi{s}}_{2,1}^{-1}(\bs{\lambda}_{2,1})}
 [\bar{\bi{s}}_1^{-1}(\bs{\lambda}_1),
       \bar{\bi{s}}_2^{-1}(\bs{\lambda}_2)]\,
\mbox{\large $|$}\partial_{\bs{\lambda}_{1,2}}\bar{\bi{s}}_{1,2}^{-1}
(\bs{\lambda}_{1,2})\mbox{\large $|$}
\end{split}
\label{eq:cftr2}
\end{equation} 
\end{linenomath}
and
\begin{linenomath} 
\begin{equation}
\begin{split}
&\eta_{I',\mathbf{\Lambda}_{1,2}|\bs{\lambda}_{2,1}}(\bi{y}_{1,2},\bs{\lambda}_{2,1})= 
\chi_{I',\mathbf{\Lambda}_{1,2}|\bs{\lambda}_{2,1}}\,\times \\
&\eta_{I,\mathbf{\Theta}_{1,2}|\bar{\bi{s}}_{2,1}^{-1}(\bs{\lambda}_{2,1})}
[\bi{s}^{-1}(\bi{y}_{1,2}),\bar{\bi{s}}_{2,1}^{-1}
 (\bs{\lambda}_{2,1})]\,
\mbox{\large $|$}\partial_{\bi{y}_{1,2}}\bi{s}^{-1}(\bi{y}_{1,2})\mbox{\large $|$} 
\end{split}
\label{eq:nftr2}
\end{equation} 
\end{linenomath}
are the transformations of the consistency and the normalization
factors that are induced by the transformations
$(\bar{\bi{s}}_1,\bar{\bi{s}}_2,\bi{s}):
(\mathbf{\Theta}_1,\mathbf{\Theta}_2,\bi{X})\longrightarrow
(\bar{\bi{s}}_1\circ\mathbf{\Theta},\bar{\bi{s}}_2\circ\mathbf{\Theta},
\bi{s}\circ\bi{X})\equiv
(\mathbf{\Lambda}_1,\mathbf{\Lambda}_2,\bi{Y})$ of the random variable
$(\mathbf{\Theta}_1,\mathbf{\Theta}_2,\bi{X})$.
\end{proper}
\hskip 0mm \textit{Proof.} \;Combining equations \eqref{eq:repar} and
\eqref{eq:condtrans3} results in
\begin{linenomath} 
\begin{equation*}
\begin{split}
&f_{I'}\!\left(\bs{\lambda}_{1,2}|\bs{\lambda}_{2,1},\bi{y}_{1,2}\right)=\\
&\frac{\zeta_{I,\mathbf{\Theta}_{1,2}|\bar{\bi{s}}_{2,1}^{-1}(\bs{\lambda}_{2,1})}
 [\bar{\bi{s}}_1^{-1}(\bs{\lambda}_1),\bar{\bi{s}}_2^{-1}(\bs{\lambda}_2)]}
 {\eta_{I,\mathbf{\Theta}_{1,2}|\bar{\bi{s}}_{2,1}^{-1}(\bs{\lambda}_{2,1})}
[\bi{s}^{-1}(\bi{y}_{1,2}),\bar{\bi{s}}_{2,1}^{-1}
 (\bs{\lambda}_{2,1})]}
 \frac{\mbox{\large $|$}\partial_{\bs{\lambda}_{1,2}}
       \bar{\bi{s}}_{1,2}^{-1}(\bs{\lambda}_{1,2})\mbox{\large $|$}}
      {\mbox{\large $|$}\partial_{\bi{y}_{1,2}}\bi{s}^{-1}(\bi{y}_{1,2})
\mbox{\large $|$}}\times\\
&f_{I'}\!\left(\bi{y}_{1,2}|\bs{\lambda}_1,\bs{\lambda}_2\right),
\end{split}
\end{equation*} 
\end{linenomath}
which, when compared to the relation
\begin{linenomath} 
\begin{equation*}
\begin{split}
&f_{I'}\!\left(\bs{\lambda}_{1,2}|\bs{\lambda}_{2,1},\bi{y}_{1,2}\right)= \\
&\frac{\zeta_{I',\mathbf{\Lambda}_{1,2}|\bs{\lambda}_{2,1}}
      (\bs{\lambda}_1,\bs{\lambda}_{2})}
     {\eta_{I',\mathbf{\Lambda}_{1,2}|\bs{\lambda}_{2,1}}
      (\bi{y}_{1,2},\bs{\lambda}_{2,1})}\, 
 f_{I'}\!\left(\bi{y}_{1,2}|\bs{\lambda}_1,\bs{\lambda}_2\right) \ ,
\end{split}
\end{equation*} 
\end{linenomath}
implied by Proposition\;\ref{theo:consistency}, yields \eqref{eq:cftr2} and
\eqref{eq:nftr2}. In the same way, \eqref{eq:cftr1} and \eqref{eq:nftr1} are
obtained if \eqref{eq:condtrans3} is replaced by \eqref{eq:condtrans2}.
%
\hfill $\Box$ 
\vskip 3mm  
For invariant families $I$ of direct probability distributions,
equations \eqref{eq:cftr1} and \eqref{eq:cftr2} reduce to
functional equations
\begin{linenomath} 
\begin{equation}
\zeta_{I,\mathbf{\Theta}}(\bs{\theta})= \chi_{I,\mathbf{\Theta}}(a)\,
\zeta_{I,\mathbf{\Theta}}[\bar{\bi{g}}_{a}^{-1}(\bs{\theta})]\,
 \mbox{\large $|$}\partial_{\bs{\theta}}\bar{\bi{g}}_{a}^{-1}(\bs{\theta})\mbox{\large $|$}
\label{eq:funeqcon1}
\end{equation} 
\end{linenomath}
and
\begin{linenomath} 
\begin{equation}
\begin{split}
&\zeta_{I,\mathbf{\Theta}_{1,2}|\bs{\theta}_{2,1}}(\bs{\theta}_1,\bs{\theta}_{2})
 = \chi_{I,\mathbf{\Theta}_{1,2}|\bs{\theta}_{2,1}}(a)\,\times \\
&\zeta_{I,\mathbf{\Theta}_{1,2}|\bar{\bi{g}}_{a,2,1}^{-1}(\bs{\theta}_{2,1})}
 [\bar{\bi{g}}_{a,1}^{-1}(\bs{\theta}_1),
       \bar{\bi{g}}_{a,2}^{-1}(\bs{\theta}_2)]\,
  \mbox{\large $|$}\partial_{\bs{\theta}_{1,2}}\bar{\bi{g}}_{a,1,2}^{-1}
  (\bs{\theta}_{1,2})\mbox{\large $|$}
\end{split}
\label{eq:funeqcon2}
\end{equation} 
\end{linenomath}
for the consistency factors $\zeta_{I,\mathbf{\Theta}}(\bs{\theta})$ and
$\zeta_{I,\mathbf{\Theta}_{1,2}|\bs{\theta}_{2,1}}(\bs{\theta}_1,\bs{\theta}_{2})$,
respectively. It should be noticed that the usual multipliers
$\chi_{I,\mathbf{\Theta}}$
and $\chi_{I,\mathbf{\Theta}_{1,2}|\bs{\theta}_{2,1}}$, up to which  
the two consistency factors are uniquely determined 
(Property\;\ref{prop:unique}), may depend on the
parameters $a$ of the transformations (on the group elements $a$), 
i.e., the consistency
factors for the parameters of invariant parametric families of
direct probability distributions are to be {\em relatively
invariant} under $\bar{\mathcal G}$.

Apart from the invariance of the consistency factors, 
invariance of a family $I$ of direct distributions under a group
${\mathcal G}$ also implies invariance of the family of 
the corresponding inverse distributions under the induced
group $\bar{\mathcal G}$.
Let, for example, $I$ be an invariant
parametric family of continuous direct probability distributions 
of a scalar random variable $X$, whose scalar parameter is denoted
by $\Theta$. Then, according to \eqref{eq:leminvar},
\begin{linenomath} 
\begin{equation}
 F_I\left(\theta|x\right) = 
\begin{cases}
  \hskip 3.2mm F_I(\bar{h}(a^{-1},\theta)|h(a^{-1},x)) 
 \hskip 3.2mm; \bar{g}'_a(\theta) > 0 \\
  1- F_I(\bar{h}(a^{-1},\theta)|h(a^{-1},x)) ; 
  \bar{g}'_a(\theta) < 0 
 \end{cases}
\hskip -4mm .
\label{eq:reslemm2}
\end{equation} 
\end{linenomath}

\subsection{Invariance under discrete groups of transformations}
\label{ss:discrete}

Under what circumstances functional equations \eqref{eq:funeqcon1} and
\eqref{eq:funeqcon2} lead to unique solutions 
$\zeta_{I,\mathbf{\Theta}}(\bs{\theta})$ and
$\zeta_{I,\mathbf{\Theta}_{1,2}|\bs{\theta}_{2,1}}(\bs{\theta}_1,\bs{\theta}_{2})$?
\begin{exam}[Parity]
\label{ex:parity}
{\rm Let a parametric family of continuous direct probability
distributions be invariant under a discrete group ${\mathcal G}$
of transformations $g_a: X\longrightarrow aX$ with $\bar{g}_a:
\Theta\longrightarrow a\Theta$ being the corresponding 
transformations from the induced group, where the underlying group $G$
consists of two elements, $a=\pm 1$.
That is, the distributions from the considered
family have (positive) {\em parity} under 
simultaneous inversions of the spaces of $X$ and $\Theta$.
By combining $\zeta_{I,\Theta}(\theta)=
\chi_{I,\Theta}(a)\,\zeta_{I,\Theta}[g_a^{-1}(\theta)]$
and $\zeta_{I,\Theta}[g_a^{-1}(\theta)]=\chi_{I,\Theta}(a)\,
\zeta_{I,\Theta}\{ g_a^{-1}[g_a^{-1}(\theta)]\}$
and setting $a=-1$ we obtain  $\zeta_{I,\Theta}(\theta)=
\chi_{I,\Theta}(-1)\,\zeta_{I,\Theta}(-\theta)$ and
$\zeta_{I,\Theta}(\theta)=
[\chi_{I,\Theta}(-1)]^2\,\zeta_{I,\Theta}(\theta)$, such that
$[\chi_{I,\Theta}(-1)]^2=1$. When inability of 
$\zeta_{I,\Theta}(\theta)$ to switch sign is invoked
(Property\;\ref{prop:sign}), this further implies
$ \zeta_{I,\Theta}(-\theta)= \zeta_{I,\Theta}(\theta)$.
That is, $\zeta_{I,\Theta}(\theta)$ must have positive 
parity under the inversion $\Theta\longrightarrow -\Theta$, 
but apart from this, it can take any
form and so in this case equation \eqref{eq:funeqcon1}
does not lead to unique solution. 
}
\end{exam}

It is not difficult to understand that
this is a common feature of all solutions based on invariance
of parametric families under \textit{discrete} groups. 
If the symmetry group is discrete, the spaces of $X$ and $\Theta$
break up in intervals, the so-called \textit{fundamental
regions} or \textit{domains} of the group
(\cite{wig}, \S\,19.1, p.\,210; \cite{jay}, \S\,10.9, p.\,332), 
with no connections in terms of group transformations within the
points of the same interval. We are then free to choose the form 
of $\zeta_{I,\Theta}(\theta)$ in one of these intervals (e.g., we can choose 
$\zeta_{I,\Theta}(\theta)$ for the positive values of $\theta$ in the above
example), hence the invariance of a family $I$ under a 
discrete group ${\mathcal G}$ alone does not lead to a unique form of
the corresponding consistency factor. 
The argument applies, for example, for all parametric families of
discrete probability distributions.
%

\subsection{Consistency factors and invariance under Lie groups}
\label{ss:lie1}

Let ${\mathcal G}=\{g_a:\mathbb{R}\longrightarrow\mathbb{R}\,;\ a\in
G\}$ be a group and $G$ be a one-dimensional Lie group. Then,
according to Proposition \ref{theo:existence},
on the subspace $\widetilde{V}_{X}\subseteq V_X$ with non-vanishing derivative
\eqref{eq:prtder1}, every ${\mathcal G}$-invariant parametric family
$I$ of continuous direct probability distributions 
is necessarily isomorphic to a location-scale family $I'$  
with the realization $\sigma=1$ of the scale parameter $\Theta_2$.
Since the fundamental domain of the group 
$\bar{\mathcal{G}}$ of translations on the real axis consists 
of a single point, the space of all possible realizations 
of a location parameter is a \textit{homogenous space for the group} 
(i.e., the space is said to be a \textit{single} 
$\bar{\mathcal{G}}$\textit{-orbit}).

The implications of Proposition\;\ref{theo:existence} may be extended 
to the subspaces $V_X-\widetilde{V}_{X}$:
\begin{theo}
\label{theo:existenceprime}
Let ${\mathcal G}=\{g_a:a\in G\}$ be a group of transformations
$g_a:\mathbb{R}\longrightarrow\mathbb{R}$
and $G$ be a one-dimensional Lie group. Suppose, in addition, that 
a parametric family $I$ 
of continuous direct probability distributions for a scalar random variable $X$  
is $\mathcal{G}$-invariant, that the action of $G$ on $\mathbb{R}$
is not identically trivial on entire $V_X$, and that
the corresponding cdf's $F_I(x|\lambda)$ are differentiable in $\lambda$. 
Then, for a realization $x\in V_X-\widetilde{V}_X\subset V_X$ 
with vanishing derivative \eqref{eq:prtder1}, the inverse probability
distribution whose cdf $F_I(\lambda|x)$ is differentiable in
$x$, cannot be assigned. (Existence of 
derivatives $\partial_xF_I(x|\lambda)$ and  $\partial_{\lambda}F_I(\lambda|x)$
is assured by Definition\;\ref{def:scpdf}.)
\end{theo}
\begin{exam}
\label{ex:zeroder}
{\rm Let
  $I_{\mu}\equiv\{P_{I,(\mu,\sigma)}:(\mu,\sigma)\in(\mu,\mathbb{R}^+)\}$ 
 be a sub-family of a continuous location-scale family $I$ that
  corresponds to the value of the location parameter $\Theta_1$ being
fixed to $\mu$. By transformation $X\longrightarrow X-\mu \equiv Y$, every cdf
$F_{I_{\mu}}\left(x|\mu,\sigma\right)$ from $I_{\mu}$ is reduced to
\begin{linenomath} 
\begin{equation*}
 F_{I'_{\mu}}\left(y|\mu,\sigma\right)= F_{I_{\mu}}\left(y+\mu|\mu,\sigma\right)
 =  \Phi\left(\frac{y}{\sigma}\right) \ ,
\end{equation*} 
\end{linenomath}
where $y\equiv x-\mu$. The probability distribution for the random
variable $Y$ thus belongs
to the family $I_{\mu}'$ that is invariant under transformations $g_a: Y
\longrightarrow a Y$ and $\bar{g}_a: \Theta_2\longrightarrow
a\Theta_2$ for all $a\in\mathbb{R}^+$. Since the derivative
$\partial_ah(a^{-1},y)\mbox{\large $|$}_{a=e}=y$ vanishes for $y=0$, 
the inverse probability distribution for the scale parameter 
$\Theta_2$ given $y=0$ (or, equivalently, given $x=\mu$) 
does not exist. 
}
\end{exam}
In order to assign an inverse probability distribution to a scalar 
parameter of a family that is invariant under a group ${\mathcal G}$ that is
underlain by a one dimensional Lie group $G$ it therefore suffices 
to determine the
consistency factor $\zeta_{I,\Theta_1}(\mu)\equiv
\zeta_{I,\Theta_1|\sigma=1}(\mu,\sigma)$, which can
subsequently be transformed, by means of \eqref{eq:cftr1}, to the 
corresponding consistency factor $\zeta_{I',\Lambda}(\lambda)$ for the original 
parameter $\Lambda$. 
A location-scale family $I_{\sigma=1}=\{
Pr_{I',(\mu,\sigma)}:(\mu,\sigma)\in(\mathbb{R},1)\}$ of continuous direct
  probability distributions with the fixed value $\sigma=1$ of the
  scale parameter is a subset of the location-scale family 
$I=\{Pr_{I,(\mu,\sigma)}:(\mu,\sigma)\in\mathbb{R}\times\mathbb{R}^+
\}$ that is
    invariant under the group ${\mathcal G}$ \eqref{eq:invtransx}. 
Given a location-scale family $I$, the functional equation 
\eqref{eq:funeqcon1} for the consistency factor 
$\zeta_{I,\Theta_1|\sigma}(\mu,\sigma)$ therefore reduces to
\begin{linenomath} 
\begin{equation}
\zeta_{I,\Theta_1|\sigma}(\mu,\sigma)=h(a_1,a_2)\,
\zeta_{I,\Theta_1|\sigma}[(\mu-a_1)/a_2,\sigma/a_2] \ ,
\label{eq:pil1}
\end{equation} 
\end{linenomath}
$\mu,a_1\in{\mathbb{R}}$ and $\sigma,a_2\in{\mathbb{R}}^+$,
where $h(a_1,a_2)\equiv \chi_{I,\Theta_1|\sigma}(a_1,a_2)/a_2$. 
\begin{lemm}
\label{lem:loc}
The solution $\zeta_{I,\Theta_1|\sigma}(\mu,\sigma)$ of equation
\eqref{eq:pil1} is a function of $\sigma$ alone, say $\Omega(\sigma)$.
\end{lemm}
Since $\zeta_{I,\Theta_1|\sigma}(\mu,\sigma)$ is uniquely determined
only up to a factor  $\chi_{I,\Theta_1|\sigma}(x_{1,2},\sigma)$
(Property\;\ref{prop:unique}), $\Omega(\sigma)$ may be, without 
loss of generality, set to unity, such that 
\begin{linenomath} 
\begin{equation}
\zeta_{I,\Theta_1|\sigma}(\mu,\sigma) = 1 \ ,
\label{eq:pimufinal}
\end{equation} 
\end{linenomath}
regardless the explicit family $I$ of direct probability
distributions, as well as the realization $\sigma$ of the scale parameter.

By using the same arguments as for
$\zeta_{I,\Theta_1|\sigma}(\mu,\sigma)$ we find that a consistency
factor $\zeta_{I,\Theta_2|\mu}(\mu,\sigma)$ is also a function of
$\sigma$ only, say
$\zeta_{I,\Theta_2}(\sigma)\equiv\zeta_{I,\Theta_2|\mu}(\mu,\sigma)$.
The inverse probability distribution for the scale parameter
$\Theta_2$, given $\Theta_1=\mu$ and $X_1=x_1=\mu$, does not exist
(Example\;\ref{ex:zeroder}), while for $x_1\gtrless\mu$ the pdf
$f_I\left(\sigma|\mu,x_1\right)$ can be expressed in terms of
$f_{I^{\pm}}\left(x_1|\mu,\sigma\right)$ (Section\;\ref{ss:parfam}):  
\begin{linenomath} 
\begin{equation*}
\begin{split}
     f_I\left(\sigma|\mu,x_1\right)&=
 \frac{\zeta_{I,\Theta_2|\mu}(\mu,\sigma)\,f_I\left(x_1|\mu,\sigma\right)}
   {\eta_{I,\Theta_2|\mu}\left(x_1,\mu\right)} \\
&=
 \frac{\zeta_{I,\Theta_2}(\sigma)\,f_{I^{\pm}}\left(x_1|\mu,\sigma\right)}
   {\eta_{I^{\pm},\Theta_2|\mu}\left(x_1,\mu\right)}
\ ,
\end{split}
\end{equation*} 
\end{linenomath}
where $\eta_{I^{\pm},\Theta_2|\mu}\left(x_1,\mu\right)\equiv
\eta_{I,\Theta_2|\mu}\left(x_1,\mu\right)/c_{\pm}$. By equation 
\eqref{eq:sc2loc}, every pdf
$f_{I^{\pm}}\left(x_1|\mu,\sigma\right)$ is reducible to 
$f_{{I^{\pm}}'}(y_1|\lambda_1,\lambda_2=1)$, such that
\begin{linenomath} 
\begin{equation*}
     f_{I'}\left(\lambda_1|\lambda_2=1,y_1\right)=
 \frac{\zeta_{I',\Lambda_1}(\lambda_1)\,
       f_{I^{\pm}{'}}\left(y_1|\lambda_1,\lambda_2=1\right)}
   {\eta_{I^{\pm}{'},\Lambda_1|\lambda_2=1}\left(y_1,\lambda_2=1\right)}
\end{equation*} 
\end{linenomath}
holds true and
$\zeta_{I',\Lambda_1}(\lambda_1)\equiv
\zeta_{I',\Lambda_1|\lambda_2=1}(\lambda_1,\lambda_2)=1$, 
where $y_1\equiv\ln{\{\pm(x_1-\mu)\} }$ and
$\lambda_1\equiv\ln{\sigma}\equiv\bar{s}(\sigma)$.
Since, according to equation \eqref{eq:cftr1},
\begin{linenomath} 
\begin{equation*}
 \zeta_{I',\Lambda_1}(\lambda_1)=\,
 \zeta_{I,\Theta_2}[\bar{s}^{-1}(\lambda_1)]\,
 \mbox{\large{$|$}}[\bar{s}^{-1}(\lambda_1)]' 
\mbox{\large$|$}
\end{equation*} 
\end{linenomath}
must also hold, 
\begin{linenomath} 
\begin{equation}
\zeta_{I,\Theta_2|\mu}(\mu,\sigma) = \sigma^{-1} 
\label{eq:pisfinal}
\end{equation} 
\end{linenomath}
is the general form of the consistency factor 
$\zeta_{I,\Theta_2|\mu}(\mu,\sigma)$,
again regardless the explicit location-scale family $I$ of direct probability
distributions and the realization $\mu$ of the location  parameter.

According to Proposition\;\ref{theo:consistency}, an inverse pdf 
$f_I\left(\mu,\sigma|x_1,x_2\right)$ for the parameters $\Theta_1$ and
$\Theta_2$ of a location-scale family $I$ must be expressible as 
\begin{linenomath} 
\begin{equation*}
     f_I\left(\mu,\sigma|x_1,x_2\right)=
 \frac{\zeta_{I,\mathbf{\Theta}}(\mu,\sigma)\,f_I\left(x_1,x_2|\mu,\sigma\right)}
   {\eta_{I,\mathbf{\Theta}}\left(x_1,x_2\right)}\ .
\end{equation*} 
\end{linenomath}
For the same reasons as $\zeta_{I,\Theta_1|\sigma}(\mu,\sigma)$ (Lemma
\ref{lem:loc}), $\zeta_{I,\mathbf{\Theta}}(\mu,\sigma)$ must also be a 
function of $\sigma$ alone, say $\Xi(\sigma)$, while the product 
rule \eqref{eq:invproduct} implies factorizability of
$ f_I\left(\mu,\sigma|x_1,x_2\right)$,
\begin{linenomath} 
\begin{equation}
    f_I\left(\mu,\sigma|x_1,x_2\right)=
    f_I(\sigma|\mu,x_1,x_2)\,f_I(\mu|x_1,x_2) ,
\label{eq:prodinverse}
\end{equation} 
\end{linenomath}
where, according to Bayes' Theorem \eqref{eq:bayes1},
\begin{linenomath} 
\begin{equation*}
\begin{split}
     f_I\left(\sigma|\mu,x_1,x_2\right)&=
 \frac{f_I\left(\sigma|\mu,x_1\right)\,f_I\left(x_2|\mu,\sigma\right)}
     {f_I\left(x_2|\mu,x_1\right)} \\
&=
 \frac{\zeta_{I,\Theta_2|\mu}(\mu,\sigma)\,f_I\left(x_1,x_2|\mu,\sigma\right)}
   {\eta_{I,\Theta_2|\mu}(x_1,\mu)\,f_I\left(x_2|\mu,x_1\right)}\ .
\end{split}
\end{equation*} 
\end{linenomath}
Hence,
\begin{linenomath} 
\begin{equation*}
 \frac{\Xi(\sigma)}
      {\eta_{I,\mathbf{\Theta}}\left(x_1,x_2\right)}=
\frac{\sigma^{-1}\,f_I\left(\mu|x_1,x_2\right)}
     {\eta_{I,\Theta_2|\mu}(x_1,\mu)\,f_I\left(x_2|\mu,x_1\right)} \ .
\end{equation*} 
\end{linenomath}
must hold, finally implying
\begin{linenomath} 
\begin{equation}
\Xi(\sigma)=\zeta_{I,\mathbf{\Theta}}(\mu,\sigma) =
\sigma^{-1} 
\ .
\label{eq:pimusfinal}
\end{equation} 
\end{linenomath}

The findings of the present subsection can thus be recapitulated as follows:
\begin{theo}
\label{theo:consfact}
The consistency factors $\zeta_{I,\Theta_1|\sigma}(\mu,\sigma)$,
$\zeta_{I,\Theta_2|\mu}(\mu,\sigma)$
and $\zeta_{I,\mathbf{\Theta}}(\mu,\sigma)$ for the parameters of
location-scale families of continuous direct probability distributions 
read $\zeta_{I,\Theta_1|\sigma}(\mu,\sigma) = 1$ and 
$\zeta_{I,\Theta_2|\mu}(\mu,\sigma) = \zeta_{I,\mathbf{\Theta}}(\mu,\sigma) 
  = \sigma^{-1}$.
\end{theo}
%

\subsection{On integrability and on uniqueness of the consistency factors}
\label{ss:unique}

It is easily verified that normalizability of pdf's
\eqref{eq:defpdflocscale} from location-scale families guarantees
also normalizability (integrability) 
of all the pdf's that were involved in the foregoing 
derivations of the consistency factors. No requirement concerning 
integrability, however, has ever been imposed to
consistency factors themselves. Moreover, it is evident that
consistency factors $\zeta_{I,\Theta_1|\sigma}(\mu,\sigma)$
\eqref{eq:pimufinal}, defined on the entire real axis, are not
integrable, implying that none of the consistency factors for
scalar parameters of parametric families 
that are invariant under the action of a one-dimensional Lie
group, is integrable. 

Let  $\zeta_{I,\mathbf{\Theta}}(\bs{\theta})$ be a non-integrable
consistency factor for a parameter $\mathbf{\Theta}$ from a 
family $I$ of continuous direct probability distributions. Suppose for
a moment that apart from the conditional pdf's
$f_I(\bi{x}|\bs{\theta})$ and $f_I(\bs{\theta}|\bi{x})$, there also 
exist the non-informative prior pdf $f_I(\bs{\theta})$
and the joint pdf $f_I(\bs{\theta},\bi{x})$. Then, there exists
an unconditional predictive pdf $f_I(\bi{x})$ (see, for
  example, \cite{sha}, \S\,4.1.1, Theorem\;4.1, p.\,194), such that
\begin{linenomath} 
\begin{equation}
 f_I(\bs{\theta}|\bi{x})=\frac{f_I(\bs{\theta})\,f_I(\bi{x}|\bs{\theta})}
                              {f_I(\bi{x})} \ .
\label{eq:bayes3}
\end{equation} 
\end{linenomath}
But apart from Bayes' Theorem \eqref{eq:bayes3},
$f_I(\bs{\theta}|\bi{x})$ is also subjected to
Proposition\;\ref{theo:consistency}, implying that $f_I(\bs{\theta})$ and
$\zeta_{I,\mathbf{\Theta}}(\bs{\theta})$ are equal up to an arbitrary
multiplication constant. Since then $f_I(\bs{\theta})$ is not
integrable, the non-informative pdf $f_I(\bs{\theta})$ does not
exist, and consequently, neither do exist $f_I(\bs{\theta},\bi{x})$
and 
the underlying probability space $(\Omega,\Sigma,P)$. The 
pdf's $f_I(\bi{x}|\bs{\theta})$ and 
$f_I(\bs{\theta}|\bi{x})$ therefore represent an extension
of the concept of the conditional probability distribution
that was introduced in Subsection\;\ref{ss:general}.

Since every consistency factor is determined only up to an arbitrary
multiplicative factor (Property \ref{prop:unique}), infinitely many 
different consistency factors for a parameter from a particular
parametric family exist. Nevertheless, unlike non-unique
non-informative prior probability distributions (recall the assertions
quoted in the introductory remarks), for a scalar parameter of
a family of direct probability distributions whose invariance is
associated to a one-dimensional Lie group, for example, 
the consistency factors are unique in that they all lead to 
the same inverse probability distribution.

\subsection{Discussion}
\label{ss:discussion}
Above, the consistency factors were deduced exclusively by presuming
existence of the inverse probability distributions and by making use 
of the invariance of the families of direct probability distributions
that is related to Lie groups. The resulting set of the families with possible
probabilistic parametric inference is limited: for example, for 
scalar random variables $X$ and scalar parameters $\Theta$ the
probabilistic parametric inference is in this way restricted to location
parameters (or to parameters that are reducible to location parameters by
one-to-one transformations). On the other hand, several principles 
were proposed for determination of the non-informative prior
probability distributions. Here, applicability of these principles
for determination of the consistency factors is investigated in order
to extend the domain of the probabilistic parametric inference. 

For example, if adapted for determination of consistency factors,
{\em Bayes' Postulate} \citep{bay}, also referred to as 
the {\em Laplace Principle of Insufficient Reason}
(\cite{lap1}, p.\,XVII), suggests that all consistency factors
should be uniform. Clearly, this is inadmissible since in general 
the constant consistency factors contradict expressions \eqref{eq:cftr1} and
\eqref{eq:cftr2} for transformations of the consistency factors under 
reparameterizations.

A sophisticated version of the Principle of Insufficient
Reason is referred to as the \textit{Principle of Maximum Entropy}. 
In our context, the {\em information entropy} (\cite{shan}, \S\,6) reads
\begin{linenomath} 
\begin{equation*}
  S\equiv -\int_{V_{\Theta}}\zeta_{I,{\Theta}}(\theta)\,
     \ln{\left\{\zeta_{I,{\Theta}}(\theta)\right\} }\,d\theta \ ,
\end{equation*} 
\end{linenomath}
while the Principle of Maximum Entropy states 
(\cite{jay}, \S\,11.3, pp.\,350) that the consistency factor 
which maximizes the entropy
represents the most honest description of what we know about
the value of the inferred parameter. For compact parameter spaces
$V_{\Theta}$ for which the above integral exists, the principle
again results in constant consistency factors 
$\zeta_{I,\Theta}(\theta)=e^{-1}$. The factors are
then flawed in the same way as the factors implied by 
Bayes' Postulate. \citet[][\S\,12.3, pp.\,374-377]{jay} 
argues that the above expression for the entropy is inappropriate
since it is not invariant under reparameterization and proposes
a {\em Kullback-Leibler divergence} (also called {\em relative entropy}) to
replace it:
\begin{linenomath} 
\begin{equation*}
  S\equiv -\int_{V_{\Theta}}\zeta_{I,{\Theta}}(\theta)\,
     \ln{\left\{\frac{\zeta_{I,{\Theta}}(\theta)}
        {m(\theta)}\right\} }\,d\theta \ , 
\end{equation*} 
\end{linenomath}
where $m(\theta)$ is the {\em reference measure
function}. Due to the unknown form of the latter, however,
maximization of the relative entropy does not
lead to unique consistency factors.

If {\em Jeffreys' general rule} is applied \cite{jef2}, the consistency
factors are determined \textit{via} the determinant of the Fisher
information matrix ${\mathcal I}_{I,\mathbf{\Theta}}(\bs{\theta})$,
$\zeta_{I,\mathbf{\Theta}}(\bs{\theta})
\propto \sqrt{\det{[ {\mathcal
    I}_{I,\mathbf{\Theta}}(\bs{\theta})]}}$, where the elements
of the matrix are given by
\begin{linenomath} 
\begin{equation*}
\begin{split}
&[{\mathcal I}_{I,\mathbf{\Theta}}(\bs{\theta})]_{i,j} \equiv \\
&\int_{\mathbb{R}^n}
 \partial_{\theta_i}\!\ln\!{\{f_I(\bi{x}|\bs{\theta})\} }\, 
 \partial_{\theta_j}\!\ln\!{\{f_I(\bi{x}|\bs{\theta})\} }\, 
 f_I(\bi{x}|\bs{\theta})\,d^n\bi{x} \ .
\end{split}
\end{equation*} 
\end{linenomath}
The obtained consistency factors satisfy requirements \eqref{eq:cftr1}
and \eqref{eq:cftr2} for transformations of the factors under
reparameterization, but are flawed in another way. Let, for
example, a probability distribution $N(\mu,\sigma)$ for a random variable $X$ 
belong to the normal (or Gaussian) family (\cite{ken}, \!\S\,5.36,
  p.\,191). Then, Jeffrey's general rule yields the 
consistency factors $\zeta_{I,\Theta_1|\sigma}(\mu,\sigma)\propto 1$,
$\zeta_{I,\Theta_2|\mu}(\mu,\sigma)\propto \sigma^{-1}$
and $\zeta_{I,\mathbf{\Theta}}(\mu,\sigma)\propto \sigma^{-2}$, such
that the resulting inverse probability distributions violate the
product rule \eqref{eq:prodinverse}.

A modification of Jeffreys' general rule by \cite{ber} called
the {\em reference prior approach} leads to violations of the same
product rule (\cite{ber}, \S\,3.3, pp.\,118-119). Also, let
$X_1,\ldots,X_n$ be independent random variables with identical
probability distribution $N(\mu,\sigma)$. Since the normal
family is a location-scale family of continuous distributions, the consistency
factor \eqref{eq:pimusfinal} yields a unique posterior pdf 
$f_{I}\left(\bar{\bi{s}}^{-1}(\lambda,\sigma)|\bi{x}\right)$ for the 
parameter $(\Theta_1,\Theta_2)$ of the distribution, whereas the
posterior pdf for $(\Lambda,\Theta_2)\equiv
\bar{\bi{s}}\circ(\Theta_1,\Theta_2)$, $\Lambda\equiv \Theta_1/\Theta_2$,
is obtained according to \eqref{eq:cftr1} (Property\;\ref{prop:trans}),
$f_{I'}(\lambda,\sigma|\bi{x})=
 f_{I}\left(\bar{\bi{s}}^{-1}(\lambda,\sigma)|\bi{x}\right)\,|
\partial_{(\lambda,\sigma)}\bar{\bi{s}}^{-1}(\lambda,\sigma)|=
 \sigma\,f_{I}\left(\bar{\bi{s}}^{-1}(\lambda,\sigma)|\bi{x}\right)$.
A unique $f_{I'}(\lambda,\sigma|\bi{x})$ further implies a unique
marginal pdf $f_{I'}(\lambda|\bi{x})$,
\begin{linenomath} 
\begin{equation*}
\begin{split}
 f_{I'}(\lambda|\bi{x})&=\int_{0}^{\infty}
 f_{I'}(\lambda,\sigma|\bi{x})\,d\sigma \\
 &\propto \exp{\left\{ -n\lambda^2/2\right\}}
 \!\int_{0}^{\infty}\hskip -2mm u^n\exp\!{\left\{\!-\frac{u^2}{2}
 +r\lambda u\right\}}\,du ,
\end{split}
\end{equation*} 
\end{linenomath}
$r\equiv(\sum x_i)/\sqrt{\sum x_i^2}$, while the reference prior
approach leads to 
\begin{linenomath} 
\begin{equation*}
\begin{split}
&f_{I'}(\lambda|\bi{x}) \propto \\
&\frac{\exp{\left\{ -n\lambda^2/2 \right\}}}
              {\sqrt{1+\lambda^2/2}}
 \int_{0}^{\infty}u^{n-1}\,\exp{\left\{-\frac{u^2}{2}+r\lambda u\right\}}\,du
\end{split}
\end{equation*} 
\end{linenomath}
(\cite{ber}, \S\,5.1, pp.\,122-123). In this way, 
since the two expressions for
$f_{I'}(\lambda|\bi{x})$ are incompatible, inconsistency of the
reference prior approach with the probabilistic parametric inference
is once more demonstrated.

Invariance theory has played an important role in the theory of
non-informative prior probability distributions (see, for example,
\cite{hart}; 
\cite{jay1} and 2003, Chapter\;12, pp.\,372-396;
\cite{daw}, Section\;2, pp.\,195-199;
\cite{vil1} and 1981; 
\cite{eat};
\cite{kas}, \S\,3.2, pp.\,1347-1348). 
Functional equations \eqref{eq:funeqcon1} and 
\eqref{eq:funeqcon2}, for instance, correspond to what has been called the
{\em Principle of Relative Invariance} \citep{hart}. Since the
relative invariance of the consistency factors is implied immediately
by the existence of the inverse probability distributions, the Principle
of Relative Invariance, when applied to consistency factors, is
redundant. Contrary to what is demonstrated above, it has also been
believed that the Principle is insufficient to determine uniquely defined priors
(consistency factors) (\cite{hart}, \S\,4, p.\,838 and \S\,10,
p.\,845; \cite{vil1}, \S\,2, p.\,454; \cite{kas},
\S\,3.2, p.\,1348).

If multipliers $\chi_{I,\mathbf{\Theta}}(a)$ and 
$\chi_{I,\mathbf{\Theta}_{1,2}|\bs{\theta}_{2,1}}(a)$ are set to
unity, equations \eqref{eq:funeqcon1} and \eqref{eq:funeqcon2}
lead to {\em inner} (or {\em form invariant}) {\em consistency factors}
(\cite{vil1}; \cite{har}, \S\,2.3, pp.\,11-12 and \S\,6.3,
pp.\,53-54). Since, however, the form invariant consistency
factors  for location-scale families,  
$\zeta_{I,\Theta_1|\sigma}(\mu,\sigma)\propto 1$,
$\zeta_{I,\Theta_2|\mu}(\mu,\sigma)\propto \sigma^{-1}$
and $\zeta_{I,\mathbf{\Theta}}(\mu,\sigma)\propto \sigma^{-2}$, lead
to a violation of the product rule \eqref{eq:prodinverse}, the {\em Principle
of Form Invariance} is inconsistent
with the probabilistic parametric inference.

When a parameter space $V_{\mathbf{\Theta}}$ of a family $I$ is identical 
to the symmetry group $G$ of the family, every realization
$\bs{\theta}$ of the parameter $\mathbf{\Theta}$ identifies both
an element of the family $I$ and an element in $G$. If, in addition,
the left action $\bar{\bi{l}}: G\times G\longrightarrow G$ coincides
with the composition of the group elements $a$ and $\bs{\theta}$,
$\bar{\bi{l}}(a,\bs{\theta})\equiv a\circ\bs{\theta}$, the form
invariant consistency factors $\zeta_{I,\mathbf{\Theta}}(\theta)$
are called the {\em left
Haar consistency factors}, where {\em left Haar} is due to the
multiplication of $\bs{\theta}$ by $a$ from the
left and due to the fact that
$\nu_{l,H}(d^{m}\bs{\theta})=\zeta_{I,\mathbf{\Theta}}(\theta)\,d^{m}\bs{\theta}$
leads to the left-invariant Haar measure 
\begin{linenomath} 
\begin{equation*}
 \nu_{l,H}(B)=\int_B \nu_{l,H}(d^{m}\bs{\theta})
\end{equation*} 
\end{linenomath}
on ${\mathcal B}^m$ (Haar, 1933), i.e., $\nu_{l,H}(a\circ B)=\nu_{l,H}(B)$
for all $a\in G$ and $B\in{\mathcal B}^m$, where 
$a\circ B\equiv\left\{a\circ b : b\in B\right\}$. 
Likewise, when $\bar{k}(a,\bs{\theta})\equiv \bs{\theta}\circ a$, 
the consistency factors that solve the functional equation
$ \zeta_{I,\mathbf{\Theta}}(\theta)=
 \zeta_{I,\mathbf{\Theta}}[\bar{k}(a^{-1},\bs{\theta})]
 \mbox{\large $|$}\partial_{\bs{\theta}}\bar{k}(a^{-1},\bs{\theta})\mbox{\large $|$} $
are called {\em right Haar consistency factors} on $G$.
When $G$ is a {\em topological group},
e.g., a Lie group, both the left and the right
Haar measures (consistency factors) exist and are unique, 
each up to a positive multiplication
constant \citep{nach}, but the two measures (consistency factors)
need not coincide.
For the location-scale families, for example, 
$\bar{k}[a,(\mu,\sigma)] =
(\mu+a_1\sigma,a_2\sigma)$ induces the right Haar consistency factor
$\zeta_{I,\mathbf{\Theta}}(\mu,\sigma)=\sigma^{-1}$ which,
in contrast to the corresponding 
left Haar factor, does not lead to the violation of the 
product rule \eqref{eq:prodinverse}. 

Several additional desirable properties are established for the right
Haar consistency factors (see
Section\;\ref{ss:calibration} below for an example). Nevertheless,
\cite{eat1,eat2,eat3}
showed that unless the symmetry groups are further restricted to, for
example, amenable groups, the probability 
distributions based on the predictive pdf's that are obtained
by applications of the right Haar consistency factors (priors)
are not generally consistent with the probability axioms. 
We cannot tell though, whether or not the 
right-invariant consistency factors
based on the restricted groups extend the collection
of families for which the probabilistic parametric inference is possible.

In summary, except possibly for the principle that identifies
consistency factors with the right Haar factors for the
underlying symmetry group $G$, all the principles discussed 
are either redundant, inconsistent with the probabilistic parametric
inference, or do not
lead to unique consistency factors.
%
\section{Interpretations of probability distributions}
\label{sec:interpretations}
\hfill \parbox{0.85\linewidth}{\small
\noindent
Every axiomatic (abstract) theory admits, as is well known, of an
unlimited number of concrete interpretations besides those from 
which it was derived. Thus we find applications in fields
of science which have no relations to the concepts 
of random event and of probability in the precise meaning of these words.
\\ \vskip -3mm
\hfill
\cite{kol}, Chapter\;1, p.\,1.
} \\
%
\subsection{Probability distributions, relative frequencies and degrees
  of belief}
\label{ss:frequencies}
So far, a mathematical theory of probabilistic parametric inference
has been discussed. In the present section, however, two
concepts of probability distributions are introduced that link the 
mathematical theory to an external world of measurable phenomena: 
the \textit{concept of relative frequencies} in repeated trials, 
and the \textit{concept of degrees of belief} in \textit{hypotheses} 
or \textit{propositions} (i.e., in statements that can be either true 
or false) concerning values of inferred parameters of parametric families.

Suppose an experiment is 
repeated under identical conditions, but the outcomes vary from one
repetition of the experiment to another. If a numerical
characteristic assigned to the outcomes of the experiment follows no describable
deterministic pattern, the experiment is called {\em random
  experiment}, the outcomes of
the experiment are called {\em random events}, while the underlying 
process of such an experiment is called {\em random process}. Let
random events be mutually independent. Then, within
the frequency interpretation of probability distributions, the direct
probability distribution for a random variable $\bi{X}$, linked
to the experiment, is assumed to coincide
with the long term distribution of relative frequencies
of particular outcomes of the experiment,
\begin{linenomath} 
\begin{equation*}
 F_{I}(\bi{x}|\bs{\theta})=\lim_{N\to\infty}\frac{N_{\bi{X}\le\bi{x}}}{N} \ ,
\end{equation*} 
\end{linenomath}
where $N$ is the total number of repetitions of the experiment and
$N_{\bi{X}\le\bi{x}}$ is the number of the repetitions with outcomes whose
numerical characteristic is less-or-equal to $\bi{x}$. Henceforth,
the frequency interpretation of direct probability distributions 
is assumed. 

Inverse probability distributions, on the other hand, are used to
express one's degrees of belief that, given a (finite) recorded sequence 
$\bi{x}_1,\bi{x}_2,\ldots$ of realizations of independent random variables
$\bi{X}_1,\bi{X}_2,\ldots$ with an identical probability distribution
from a parametric family $I$, the so-called {\em true value} of the
parameter $\mathbf{\Theta}$ of the family (i.e., the value of the
parameter that uniquely determines the true limiting frequency distribution
of the realizations) lies within a certain region of the parameter space.
Several strong arguments exist for inverse probability distributions
being the ideal for parametric inferences, like, for example, the
so-called Dutch Book Theorem, emerging from the work of 
\citet[][Chapter\;VII, pp.\,156-198]{ram}, 
\cite{def1,def2}, \cite{shi} and \cite{kem},
and Cox's Theorem \citep{cox}. For a concise review of the two Theorems
see, for example, \cite{par}, Chapter\;3, pp.\,19-33.

While being identical objects from a mathematical perspective, the
direct and the inverse probability distributions obviously have different 
interpretations. Contrary to the distribution of realizations of random
variables $\bi{X}_i$, in most situation the realization of a parameter
-- the inferred true value of $\mathbf{\Theta}$ -- 
is unknown but fixed. Several authors overlooked this important
difference between the frequency distributions and the distributions
of someone's beliefs (see, for example, 
\cite{leh}, \S\,1.6, p.\,14; \cite{sha}, \S\,7.1.3, p.\,431; 
\cite{case}, \S\,7.2.3, p.\,324 and \S\,9.2.4, pp.\,435-436;
\cite{har}, \S\,2.5, p.\,18). It should be noticed,
however, that the developed theory of probabilistic parametric 
inference still provides verifiable predictions in terms of
relative frequencies of confidence intervals, covering the true value
of the parameter (see Section\;\ref{ss:calibration}, below).
The theory is then both operational and objective.
%

\subsection{Calibration}
\label{ss:calibration}

\begin{define}[Confidence intervals]
\label{def:confint}
Let $f_I(\theta|x)$ be a pdf of a probability distribution for
a scalar parameter $\Theta$, $V_{\Theta}=
(\theta_a,\theta_b)$, 
given realization $x$ of a scalar random variable $X$ from a parametric 
family $I$. A {\em confidence interval}
$(\theta_1(x),\theta_2(x)) \subseteq V_{\Theta}$ 
is defined via the system of equations
\begin{linenomath} 
\begin{equation*}
  Pr_I(\theta_a,\theta_{1})|x) =
 \!\int_{\theta_a}^{\theta_{1}}\!\! f_I(\theta|x) 
 \,d\theta = \alpha
\end{equation*} 
\end{linenomath}
and
\begin{linenomath} 
\begin{equation*}
Pr_I((\theta_{1},\theta_{2})|x) =
 \!\int_{\theta_{1}}^{\theta_{2}}\!\!
 f_I(\theta|\bi{x})\,d\theta = 
 \delta \ ,
\end{equation*} 
\end{linenomath}
where $\delta\in[0,1]$ and $\alpha\in[0,1-\delta]$. The number
$\delta$ is called the {\em probability content} of the interval.
\end{define}
Higher dimensional {\em confidence regions}, 
e.g., $m$-dimensional {\em confidence
rectangles} ($m\ge 2$), for vector-parameters are defined in a similar way.
\begin{define}[Calibration]
\label{def:calibration}
Let ${\bi{x}_1,\ldots,\bi{x}_n}$ be a set of realizations of independent
continuous random variables $\bi{X}_1,\ldots,\bi{X}_n$ from a parametric family
$I$ of direct probability distributions. The inverse probability
distributions, assigned to the inferred parameter $\mathbf{\Theta}$ of the
family $I$, given realizations $\bi{x}_i$, are called {\em calibrated}
if, in the limit $n\to\infty$, the {\em coverage} of the corresponding
confidence regions (i.e., the relative frequency of the
regions that cover the true values of the
inferred parameter) coincides with the
probability content $\delta$ of the region.
\end{define}
Calibration of probability distributions
for inferences about location and scale parameters
is guaranteed by the fact that the consistency factors
$\zeta_{I,\theta_1|\sigma}(\mu,\sigma)$,
$\zeta_{I,\theta_2|\mu}(\mu,\sigma)$ and
$\zeta_{I,\mathbf{\theta}}(\mu,\sigma)$, determined in
Subsection\;\ref{ss:lie1}, coincide with the right Haar factors for
the group $\mathbb{R}$ for summations, for the group $\mathbb{R}^+$ for
multiplications, and for the group $\mathbb{R}\times\mathbb{R}^+$ 
for operations \eqref{eq:operlocscale}, respectively
\cite{stei,chang}.
That is to say, the resulting confidence regions
coincide with the so-called {\em classical
confidence regions}, first propounded by \cite{ney}.
It should be noticed that this holds true 
even if the true value of the inferred parameter 
arbitrarily varies from realization of one random variable to another.

It can further be shown that the consistency factors for location
and scale parameters, determined in
Subsection\;\ref{ss:lie1}, provide for a simple frequency
interpretation of the predictive distributions.

To relate probabilistic parametric inference to another concept -- that of
 the {\em fiducial inference} -- let $F_I(x|\lambda)$ be a cdf 
for a continuous one-dimensional random variable $X$ that is
either strictly increasing or strictly decreasing in a scalar
parameter $\lambda$. Then, 
a sufficient condition for an inverse probability distribution to be 
calibrated -- the so-called
\textit{fiducial condition} by \citet[][\S\,3.6, p.\;70]{fis}  -- reads:
\begin{linenomath} 
\begin{equation}
 f_I(\lambda|x) = \mbox{\large $|$}\partial_{\lambda}
 F_I(x|\lambda)\mbox{\large $|$} \ .
\label{eq:answer}
\end{equation} 
\end{linenomath}
Observe that for the inverse
pdf's, assigned to location and scale parameters by using the consistency
factors \eqref{eq:pimufinal} and \eqref{eq:pisfinal}, the condition
\eqref{eq:answer} is satisfied. Also, it is easily shown that 
congruence with the fiducial condition is preserved under 
updating that is made in accordance with Bayes' Theorem. 

Conformity with the fiducial condition \eqref{eq:answer} is invariant
under one-to-one transformations $Y\equiv s\circ X$ and 
$\Theta\equiv\bar{s}\circ\Lambda$ with non-vanishing derivatives
$\bar{s}'(\theta)$:
\begin{linenomath} 
\begin{equation*}
\begin{split}
 f_{I'}(\theta|y)&=f_I(\bar{s}^{-1}(\theta)|s^{-1}(y))\,
                 \mbox{\large $|$} [\bar{s}^{-1}(\theta)]'\mbox{\large
                   $|$} \\
&=
 \mbox{\large $|$}\partial_{\bar{s}^{-1}(\theta)}
 F_I(\bar{s}^{-1}(\theta)|s^{-1}(y))\,
 [\bar{s}^{-1}(\theta)]'\mbox{\large $|$}
\end{split}
\end{equation*} 
\end{linenomath}
and therefore
\begin{linenomath} 
\begin{equation*}
 f_{I'}(\theta|y)= \mbox{\large $|$}\partial_{\theta} F_{I'}(\theta|y)
\mbox{\large $|$} \ ,
\end{equation*} 
\end{linenomath}
where the last equality is due to equation
\begin{linenomath} 
\begin{equation*}
 F_{I'}(\theta|y) =
 \begin{cases}
 \hskip 3.2mm F_I(\bar{s}^{-1}(\theta)|s^{-1}(y))
 \hskip 3.2mm \ ; \ \bar{s}'(\theta) > 0 \\
 1- F_I(\bar{s}^{-1}(\theta)|s^{-1}(y)) \ ; 
 \ \bar{s}'(\theta) < 0
 \end{cases}
\end{equation*} 
\end{linenomath}
that follows immediately from the definition of the inverse cdf's 
and from equation \eqref{eq:leminvar}.
In addition, by combining equation \eqref{eq:bayprime2} from
Proposition\;\ref{theo:consistency} with the above fiducial condition we obtain:
\begin{linenomath} 
\begin{equation}
 \zeta_{I,\Lambda}(\lambda)\,\partial_{x}F_I(x|\lambda)\pm
 \eta_{I,\Lambda}(x)\,\partial_{\lambda}F_I(x|\lambda) = 0 \ ,
\label{eq:F1F2}
\end{equation} 
\end{linenomath}
where the upper (lower) sign stands for cdf's
which are strictly decreasing (increasing) in $\lambda$. By defining
$H(x,\lambda) \equiv s(x) \mp \bar{s}(\lambda)$, 
with $s(x)$ and $\bar{s}(\lambda)$ being related to
$\zeta_{I,\Lambda}(\lambda)$ and $\eta_{I,\Lambda}(x)$ as
$s'(x)\equiv\eta_{I,\Lambda}(x)$ and 
$\bar{s}'(\lambda)\equiv\zeta_{I,\Lambda}(\lambda)$,
functional equation \eqref{eq:F1F2} can be reduced to \eqref{eq:qpde}.
Recall that the most general solution $F_I(x|\lambda)$ of equation 
\eqref{eq:qpde} implies existence of a cdf $F_{I'}(y|\mu)$ for 
$Y\equiv s\circ X$ from a location-scale family $I'$ with 
$\mu\equiv \pm \bar{s}(\lambda)$ 
being a realization of the location parameter $\Theta_1\equiv
\bar{s}\circ\Lambda$, whereas
the scale parameter $\Theta_2$ of the family $I'$ is set to 1.
That is, the fiducial condition \eqref{eq:answer} and the requirement
\eqref{eq:bayprime2} of Proposition\;\ref{theo:consistency} combined imply
reducibility of an inferred parameter to a location parameter. 
(\cite{lin} obtained the same result
by combining the calibration condition \eqref{eq:answer} and
Bayes' Theorem \eqref{eq:bayes3}.)
For scalar parameters, the
consistency factors that were deduced on the basis of invariance of 
parametric families under the action of one-dimensional Lie groups 
are therefore the only consistency factors for which the resulting inverse
probability distributions satisfy the fiducial condition \eqref{eq:answer}.
%
%
\section{Conclusions}
\label{sec:conclusions}
For scalar parameters, invariance of a parametric family of direct
probability distributions under the action of a one-dimensional Lie group
leads to unique inverse probability distributions. The concept of
invariance is equivalent to the concept of fiducial distributions,
combined with implications of Proposition\;\ref{theo:consistency}: both
concepts lead to identical inverse distributions and are applicable
under the same conditions. When this is observed, the original idea of 
\cite{bay} and \cite{lap1} of embedding parametric 
inference in the framework of probability theory becomes perfectly 
compatible with the concept of the classical
confidence intervals \citep{ney} and with the concept of the
fiducial distributions \citep{fis2}.
Therefore, provided that adherents of the Bayesian schools of parametric
inference are willing to give up the notion of non-informative prior
probability distributions, while at the same time adherents of the
frequentist schools are willing to adopt a broader concept of random
variable that leads to existence of inverse probability distributions,  
a reconciliations between different paradigms can be reached,
probably the same kind of reconciliation that Kendall \cite{ken3} had in mind
when he wrote:
      ``Neither party can avoid ideas of the other in order to 
        set up and justify a comprehensive theory.'' 
%
\appendix

\section{Proofs of Propositions and Lemmata
}\label{app}
%
%
\subsection{ Proof of Proposition \ref{theo:condpdf}} 
The left-hand side of \eqref{eq:defcondprdis} can be rewritten as
\begin{linenomath} 
\begin{equation}
\begin{split}
 \widetilde{\nu}_{\mathbf{1}_{\bi{Y}^{-1}(U)}}(S) &=
 \int_{\bi{Z}^{-1}(S)}\bi{1}_{\bi{Y}^{-1}(U)}(\omega)\,dP(\omega) \\
 & =\int_{\Omega}
 \bi{1}_{\bi{Y}^{-1}(U)}(\omega)\,\mathbf{1}_{\bi{Z}^{-1}(S)}(\omega)\,dP(\omega) \\
 &=\int_{\mathbb{R}^{n}\times\mathbb{R}^{m}}
 \mathbf{1}_{U}(\bi{y})\,\mathbf{1}_{S}(\bi{z})\,dPr_{\bi{X}}(\bi{y},\bi{z}) \\ 
 &=\int_{U\times S}
  f_{\bi{X}}(\bi{y},\bi{z})\,d^n\bi{y}\,d^m\bi{z} \\
 &=\int_{U\times S}
  f_{\bi{X}}(\bi{z})\,\frac{f_{\bi{X}}(\bi{y},\bi{z})}{f_{\bi{X}}(\bi{z})}\,
  d^n\bi{y}\,d^m\bi{z} \\
 &= \int_{S}h(\bi{z})\,f_{\bi{X}}(\bi{z})\,d^m{\bi{z}},
\end{split}
\label{eq:long}
\end{equation} 
\end{linenomath}
$U\in{\mathcal B}^n$ and $S\in\widetilde{\mathcal B}^m$, where 
$\widetilde{\mathcal B}^m$ is a restriction of 
${\mathcal B}^m$ to $V_{\bi{Z}}$ while
\begin{linenomath} 
\begin{equation*}
h(\bi{z})\equiv \int_{U} \frac{f_{\bi{X}}(\bi{y},\bi{z})}
              {f_{\bi{X}}(\bi{z})}\,d^n\bi{y} \ .
\end{equation*} 
\end{linenomath}
In \eqref{eq:long}, the first equality follows from the definition of 
$\widetilde{\nu}_{\mathbf{1}_{\bi{Y}^{-1}(U)}}(S)$ 
(Definition\;\ref{def:condprdis}),
the third equality follows from the change of variables Theorem
(\cite{dud}, \S\,4.1, p.\,92), while the last equality follows
from Fubini's Theorem 
(\cite{bar}, Chapter\;10, pp.\,119-120). Inserting
\eqref{eq:defcondpdf} into the right-hand side
of \eqref{eq:defcondprdis} yields, on the other hand,
\begin{linenomath} 
\begin{equation*}
\int_S\left[\int_Uf_{\bi{X}}(\bi{y}|\bi{z})\,d^n\bi{y}\right]
           \!f_{\bi{X}}(\bi{z})\,d^m{\bi{z}} = 
\int_S k(\bi{z})\, f_{\bi{X}}(\bi{z})\,d^m{\bi{z}} . 
\end{equation*} 
\end{linenomath}
Let $S_{1,2}\equiv\{\bi{z}:h(\bi{z})\gtrless k(\bi{z})\}$.
Then, the equality of $h(\bi{z})$ and $k(\bi{z})$ 
$Pr_{\bi{X}}^{\bi{Z}}$-almost everywhere on $V_{\bi{Z}}$ follows
immediately from Fatou's Lemma (see, for example, \cite{bar}, 
    Chapter\;4, Corollary\;4.10 of Fatou's Lemma,
    pp.\,34-35), while the equality of 
    $f_{\bi{X}}(\bi{y},\bi{z})/f_{\bi{X}}(\bi{z})$ and 
    $f_{\bi{X}}(\bi{y}|\bi{z})$ $\nu_L$-almost everywhere on
    $\mathbb{R}^n$ is obtained in an analogous way. 
%
%
\subsection{Proof of Lemma \ref{lem:trivial}}
Let $r: \mathbb{R}\times G\longrightarrow\mathbb{R}$, $r(x,a)\equiv
l(a^{-1},x)$, be the {\em right action} of $G$ on $\mathbb{R}$. Then, 
$r(x,a\circ b) = r[r(a,x),b]$ holds true for all $a,b\in G$
and for all $x\in\mathbb{R}$. A differentiation of $r(x,a\circ b)$
with respect to $a$ thus yields
\begin{linenomath} 
\begin{equation*}
 \partial_{a{\scriptscriptstyle \circ} b}r(x,a\circ b)\,
 \partial_{a}(a\circ b)=
 \partial_{r(x,a)}r[r(x,a),b]\,
 \partial_{a}r(x,a) ,
\end{equation*} 
\end{linenomath}
which for  $b=a^{-1}$ reduces to
\begin{linenomath} 
\begin{equation*}
\begin{split}
 &\partial_{c}r(x,c)\mbox{\large $|$}_{c=e}
 \partial_{a}(a\circ b)\mbox{\large $|$}_{b=a^{-1}}
 = \\
 &\partial_{r(x,a)}r[r(x,a),b]\mbox{\large $|$}_{b=a^{-1}}
 \,\partial_{a}r(x,a) \ ,
\end{split}
\end{equation*} 
\end{linenomath}
$c\equiv a\circ b$.
The left-hand side of the above equation is zero due to the premise of
the Lemma,
\begin{linenomath} 
\begin{equation*}
  \partial_{c}r(x,c)\mbox{\large $|$}_{c=e} \equiv 
  \partial_{c}l(c^{-1},x)\mbox{\large $|$}_{c=e} = 0 \ .
\end{equation*} 
\end{linenomath}
On the right-hand side, however, the first term,
\begin{linenomath} 
\begin{equation*}
 \partial_{r(x,a)}r[r(x,a),b]\mbox{\large $|$}_{b=a^{-1}}= 
 \partial_{y}l(a^{-1},y) \equiv
 \partial_{y}g_{a^{-1}}(y)
\end{equation*} 
\end{linenomath}
is non-vanishing for all admissible values of the index $a$ and for
all real $y\equiv g_a(x)$ since differentiability of
$l(a,x)$ with respect to $x$ for every $a$ is assumed.
Then,
$\partial_{a}r(x,a)\equiv\partial_{a}l(a^{-1},x) = 0$
is implied for all permissible $a$, i.e., 
$g^{-1}_a(x)$ is permitted to depend on $x$ only, say $g_a(x)\equiv h(x)$.
When $g_e(x)=x$ is invoked, this further means $h(x)=x$ and the Lemma
is proved.
%
\subsection{Proof of Lemma \ref{lem:extheo4}}
Suppose there exists a realization
$\lambda_0$ of ${\Lambda}$ for which the partial derivative
\eqref{eq:prtder2} vanishes. Since the
family $I$ of direct distributions is invariant under $\mathcal{G}$, equation
\eqref{eq:reslemm1} applies which, when
differentiated with respect to $a$ and set afterwards $a=e$, yields 
\begin{linenomath} 
\begin{equation*}
\begin{split}
 \partial_{x}F_I(x|\lambda)\,
 \partial_{a}l(a^{-1},x)\mbox{\large $|$}_{a=e} \ \\ 
 = -
 \partial_{\lambda}F_I(x|\lambda)\,
 \partial_{a}\bar{l}(a^{-1},\lambda)\mbox{\large $|$}_{a=e} .
\end{split}
\end{equation*} 
\end{linenomath}
The second term on right-hand side of the above equation vanishes
for $\lambda=\lambda_0$, which implies
\begin{linenomath} 
\begin{equation*}
 \partial_{a}l(a^{-1},x)\mbox{\large $|$}_{a=e} = 0 
 \ \ \ \ ; \ \ \ \ \forall\,x\in V_{X} \ .
\end{equation*} 
\end{linenomath}
This means, according to Lemma\;\ref{lem:trivial}, that all
transformations $g_a\in{\mathcal{G}}$ are trivial for all 
$x\in V_X(\lambda)$, which is in direct contradiction with the 
initial premises, so that the proof is completed.
%
%
\subsection{Proof of Lemma \ref{lem:qpde}} 
It is easily shown that every cdf $F_I(x|\lambda)$ of the form 
\eqref{eq:lemqpde} solves \eqref{eq:qpde}. In order to demonstrate
that the cdf's of the form \eqref{eq:lemqpde} are also the only solutions of
\eqref{eq:qpde}, suppose for a moment that $F_I(x|\lambda)$ can be
written in terms of two independent variables,
$H(x,\lambda)$ \eqref{eq:defG} and 
$K(x,\lambda)\equiv s(x)+\bar{s}(\lambda)$,
\begin{linenomath} 
\begin{equation}
 F_I(x|\lambda)=\Phi[H(x,\lambda),K(x,\lambda)] \ ,
\label{eq:deffgh}
\end{equation} 
\end{linenomath}
where the functions $s(x)$ and $\bar{s}(\lambda)$ are defined \textit{via}
\eqref{eq:defhk1} and \eqref{eq:defhk2}. Inserting \eqref{eq:deffgh} into \eqref{eq:qpde} yields
\begin{linenomath} 
\begin{equation*}
\begin{split}
 &\partial_{K}\Phi(H,K)\,
 \left[\partial_{x}H\,\partial_{\lambda}K
 -\partial_{\lambda}H\,\partial_{x}K\right]
 = \\
&2\,s'(x)\,\bar{s}'(\lambda)\,\partial_{K}\Phi(H,K) = \\ 
& 0 \ .
\end{split}
\end{equation*} 
\end{linenomath}
Therefore, for $s'(x),\bar{s}'(\lambda)\neq 0$, 
$\partial_K\Phi(H,K)$ must vanish
identically, such that the form \eqref{eq:lemqpde} of $F_I(x|\lambda)$
is implied. If, on the other hand, any of $s'(x)$ and
$\bar{s}'(\lambda)$ vanishes, $H(x,\lambda)$ and $K(x,\lambda)$ cease to
be independent, i.e., $K(x,\lambda)=K\left[H(x,\lambda)\right]$, such that
\eqref{eq:lemqpde} again holds true, but since in this case 
$F_I(x|\lambda)$ is either a function
of $x$ alone, a function of $\lambda$ alone, or a constant, 
such a solution is inadmissible for a cdf from a parametric family.
%
%
\subsection{Proof of Proposition 
\ref{theo:consistency}}
According to the premises of the Proposition, a positive
$f_I(\bs{\theta}_1|\bs{\theta}_2,\bi{x}_1,\bi{x}_2)$
exists and can be decomposed according to \eqref{eq:bayes1}.
Let ${\bs{\theta}_1}'\in 
\widetilde{V}_{\mathbf{\Theta}_1}(\bi{x}_{1,2},\bs{\theta}_2)$ 
be another re\-ali\-zation 
of $\mathbf{\Theta}_1$ fulfilling
the conditions of the Proposition, such that
\begin{linenomath} 
\begin{equation*}
\begin{split}
 f_I(\bs{\theta}_{1}'|\bs{\theta}_2,\bi{x}_1,\bi{x}_2)&=
 \frac{f_I({\bs{\theta}_1}'|\bs{\theta}_2,\bi{x}_2)\,
 f_I(\bi{x}_1|{\bs{\theta}_1}',\bs{\theta}_2)}
      {f_I(\bi{x}_1|\bs{\theta}_2,\bi{x}_2)} \\
 &=
 \frac{f_I({\bs{\theta}_1}'|\bs{\theta}_2,\bi{x}_1)\,
       f_I(\bi{x}_2|\bs{\theta}_1',\bs{\theta}_2)}
      {f_I(\bi{x}_2|\bs{\theta}_2,\bi{x}_1)} 
\end{split}
\end{equation*} 
\end{linenomath}
is also positive. Dividing the above equation with \eqref{eq:bayes1} yields
\begin{linenomath} 
\begin{equation}
 \frac{\kappa(\bi{x}_1,{\bs{\theta}_1}',\bs{\theta}_2)}
      {\kappa(\bi{x}_1,\bs{\theta}_1,\bs{\theta}_2)} =
 \frac{\kappa(\bi{x}_2,{\bs{\theta}_1}',\bs{\theta}_2)}
      {\kappa(\bi{x}_2,\bs{\theta}_1,\bs{\theta}_2)} \ ,
\label{eq:difft1}
\end{equation} 
\end{linenomath}
$\kappa(\bi{x}_{1,2},
\bs{\theta}_{\hskip -0.3mm 1}\left.\right.\!\!\!\!\!^({'}\left.\right.\!\!\!^) ,
\bs{\theta}_2) \equiv 
f_I(\bs{\theta}_{\hskip -0.3mm 1}
\left.\right.\!\!\!\!\!^({'}\left.\right.\!\!\!^)
|\bs{\theta}_2,\bi{x}_{1,2})/
f_I(\bi{x}_{1,2}|\bs{\theta}_{\hskip -0.3mm 1}
\left.\right.\!\!\!\!\!^({'}\left.\right.\!\!\!^),\bs{\theta}_2)$.
Clearly, in order to ensure equality in \eqref{eq:difft1} for all
$\bi{x}_1$ and $\bi{x}_2$ for which 
$f_I(\bi{x}_{1,2}|\bs{\theta}_{\hskip -0.3mm 1}
\left.\right.\!\!\!\!\!^({'}\left.\right.\!\!\!^),\bs{\theta}_2)> 0$, 
the left-hand and the right-hand 
side of the equation must be independent of $\bi{x}_1$ and $\bi{x}_2$, 
but may depend on $\bs{\theta}_1$, ${\bs{\theta}_1}'$, and
$\bs{\theta}_2$:
\begin{linenomath} 
\begin{equation*}
 q(\bs{\theta}_1,{\bs{\theta}_1}',\bs{\theta}_2)\equiv 
 \frac{\kappa(\bi{x}_{1,2},{\bs{\theta}_1}',\bs{\theta}_2)}
      {\kappa(\bi{x}_{1,2},\bs{\theta}_1,\bs{\theta}_2)}
 \ .
\end{equation*} 
\end{linenomath}
The function $q(\bs{\theta}_1,{\bs{\theta}_1}',\bs{\theta}_2)$ is
factorizable,
\begin{linenomath} 
\begin{equation*}
 \frac{\zeta_{I,\mathbf{\Theta}_1|\bs{\theta}_2}({\bs{\theta}_1}',\bs{\theta}_2)}
      {\zeta_{I,\mathbf{\Theta}_1|\bs{\theta}_2}(\bs{\theta}_1,\bs{\theta}_2)}
\equiv q(\bs{\theta}_1,{\bs{\theta}_1}',\bs{\theta}_2) \ ,
\end{equation*} 
\end{linenomath}
such that
\begin{linenomath} 
\begin{equation*}
\begin{split}
 \eta_{I,\mathbf{\Theta}_1|\bs{\theta}_2}(\bi{x}_{1,2},\bs{\theta}_2)
 &\equiv 
 \frac{\zeta_{I,\mathbf{\Theta}_1|\bs{\theta}_2}({\bs{\theta}_1}',\bs{\theta}_2)}
      {\kappa(\bi{x}_{1,2},{\bs{\theta}_1}',\bs{\theta}_2)} \\
 &= \frac
 {\zeta_{I,\mathbf{\Theta}_1|\bs{\theta}_2}(\bs{\theta}_1,\bs{\theta}_2)}
 {\kappa(\bi{x}_{1,2},{\bs{\theta}_1},\bs{\theta}_2)}
 \ ,
\end{split}
\end{equation*} 
\end{linenomath}
which proves equation \eqref{eq:bayprime1}, while equation
\eqref{eq:bayprime2} is proved in a similar way by invoking
\eqref{eq:bayes2} instead of \eqref{eq:bayes1}.
%
%
\subsection{Proof of Proposition
   \ref{theo:existenceprime}}
Suppose for a moment that a pdf for
$\theta$, $f_I(\theta|x)$, can be assigned to 
$\theta\in V_{\Theta}$ based on  $x\in V_X-\widetilde{V}_{X}$ 
for which partial derivative \eqref{eq:prtder1} vanishes. 
Since the family $I$ of direct probability
distributions is ${\mathcal G}$-invariant, the distributions
assigned to $\Theta$ are invariant under the induced group 
$\bar{\mathcal{G}}$ such that equation \eqref{eq:reslemm2} applies. 
When differentiated with respect to $a$ and set afterwards $a=e$, 
\eqref{eq:reslemm2} further implies
\begin{linenomath} 
\begin{equation*}
 \partial_{x} F_I(\theta|x)\,\partial_{a}l(a^{-1},x)\mbox{\large
   $|$}_{a=e}
\! = \! -
 \partial_{\theta}F_I(\theta|x)
 \,\partial_{a}\bar{l}(a^{-1},\theta)\mbox{\large $|$}_{a=e} 
\end{equation*} 
\end{linenomath}
for all $\theta\in V_\Theta$. The left-hand side of the above equation
vanishes due to the premises, adopted at the beginning
of the proof. Since, by Lemma\;\ref{lem:extheo4}, the second term on the 
right-hand side does not vanish anywhere on $V_\Theta$, 
$\partial_{\theta}F_I(\theta|x)=f_I(\theta|x)$ must vanish
for all $\theta\in V_{\Theta}$, which is incompatible with the
normalization requirement \eqref{eq:norpdf}. 
Therefore, the assumed existence of 
$f_I(\theta|x)$, based on $x$ with vanishing derivative 
\eqref{eq:prtder1}, inevitably leads to inconsistencies and is
thus ruled out.
%
\subsection{Proof of Lemma \ref{lem:loc}}
Equation \eqref{eq:pil1} holds true for all $\mu,a_1\in\mathbb{R}$
and for all $\sigma,a_2\in\mathbb{R}^+$. For
$a_1=\mu$ and $a_2=\sigma$ we obtain
$h(\mu,\sigma)=\zeta_{I,\Theta_1|\sigma}(\mu,\sigma)/
\zeta_{I,\Theta_1|\sigma}(0,1)$,
while setting $a_1=\mu$ and $a_2=1$ reveals factorizability of 
$\zeta_{I,\Theta_1|\sigma}(\mu,\sigma)$:
\begin{linenomath} 
\begin{equation}
 \zeta_{I,\Theta_1|\sigma}(\mu,\sigma)=\frac{\zeta_{I,\Theta_1|\sigma}(\mu,1)\,
 \zeta_{I,\Theta_1|\sigma}(0,\sigma)}{\zeta_{I,\Theta_1|\sigma}(0,1)} \ .
\label{eq:factor}
\end{equation} 
\end{linenomath}
By taking these findings into account, equation \eqref{eq:pil1}
reduces to
\begin{linenomath} 
\begin{equation*}
\begin{split}
 &\zeta_{I,\Theta_1|\sigma}(\mu,1)\,\zeta_{I,\Theta_1|\sigma}(0,\sigma)\,
 [\zeta_{I,\Theta_1|\sigma}(0,1)]^2= \\
 &\zeta_{I,\Theta_1|\sigma}(a_1,1)\,\zeta_{I,\Theta_1|\sigma}(0,a_2)\,\times \\
 &\zeta_{I,\Theta_1|\sigma}[(\mu-a_1)/a_2,1]\,
 \zeta_{I,\Theta_1|\sigma}(0,\sigma/a_2) \ ,
\end{split}
\end{equation*} 
\end{linenomath}
which for $a_1=0$ and $a_2=\sigma$ yields
%
$ \zeta_{I,\Theta_1|\sigma}(\mu,1)= 
  \zeta_{I,\Theta_1|\sigma}(\mu/\sigma,1)$.
%
Hence, $\zeta_{I,\Theta_1|\sigma}(\mu,1)$ must be a constant, such
that, according to \eqref{eq:factor},
$\zeta_{I,\Theta_1|\sigma}(\mu,\sigma)$ is a function of $\sigma$ alone.
%
%
%
\bibliography{paper}


\end{document}